\documentclass[11 pt,leqno]{article}
\title{
Extensions of Johnson's and Morita's homomorphisms that map to finitely generated abelian groups
}
\date{October 25, 2009}
\author{Matthew B. Day}
\usepackage{slashbox}
\usepackage{amsmath}
\usepackage{amsfonts}
\usepackage{amssymb}
\usepackage{amsthm}
\usepackage{latexsym}
\usepackage{epsfig}
\usepackage{xypic}
\usepackage{verbatim}
\usepackage{mathrsfs}

\theoremstyle{plain} \newtheorem{theorem}{Theorem}[section]
\theoremstyle{plain} \newtheorem{proposition}[theorem]{Proposition}
\theoremstyle{plain} \newtheorem{lemma}[theorem]{Lemma}
\theoremstyle{plain} \newtheorem{keylemma}[theorem]{Key Lemma}
\theoremstyle{plain} 
\theoremstyle{plain} \newtheorem{corollary}[theorem]{Corollary}
\theoremstyle{plain} 
\theoremstyle{plain} 
\theoremstyle{plain} 
\theoremstyle{plain} 
\theoremstyle{plain} 
\theoremstyle{remark} \newtheorem{example}[theorem]{Example}
\theoremstyle{definition} 
\theoremstyle{definition} \newtheorem{definition}[theorem]{Definition}
\theoremstyle{definition} \newtheorem{theorem-definition}[theorem]{Theorem-Definition}
\theoremstyle{definition} 
\theoremstyle{plain} 
\theoremstyle{plain} 
\theoremstyle{plain} 
\theoremstyle{plain} 

\numberwithin{equation}{section}

\newcommand\Hom{\mathrm{Hom}}
\newcommand\Diff{\mathrm{Diff}}

\newcommand\Aut{\mathrm{Aut}\,}

\newcommand\IA{\mathrm{IA}}

\newcommand\R{\mathbb{R}}

\newcommand\Z{\mathbb{Z}}
\newcommand\Q{\mathbb{Q}}

\newcommand\into\hookrightarrow

\def\co{\colon\thinspace}

\newcommand\lv[1]{\mathcal{L}_{#1}}
\newcommand\lag[1]{\mathfrak{g}_{#1}}
\newcommand\lal[1]{\mathfrak{l}_{#1}}

\newcommand\bdI[1]{\mathcal{I}(#1)}
\DeclareMathOperator{\bdM}{Mod}
\newcommand\latt{\Lambda}

\newcommand\Aker[2]{\mathcal{A}_{#1}(#2)}

\newcommand{\johom}[1]{\tau_{#1}}
\newcommand{\mohom}[1]{\tilde\tau_{#1}}
\newcommand{\fjh}[1]{\delta_{#1}}
\newcommand{\mathid}{\mathrm{id}}
\newcommand{\lcs}[2]{{#1}_{(#2)}}
\newcommand{\trunc}[1]{\Gamma_{#1}}
\newcommand{\lprop}{L}
\newcommand\jhob{\xi}
\newcommand{\vind}{\hat v}
\newcommand{\lmult}[1]{L_{#1}}
\newcommand{\psmo}{\epsilon_k}

\newcommand\vct[1]{\mathbf{#1}}
\DeclareMathOperator{\Range}{Range}

\begin{document}
\maketitle

\begin{abstract}
We extend each higher Johnson homomorphism to a crossed homomorphism from the automorphism group of a finite-rank free group to a finite-rank abelian group.
We also extend each Morita homomorphism to a crossed homomorphism from the mapping class group of once-bounded surface to a finite-rank abelian group.
This improves on the author's previous results~\cite{day}.
To prove the first result, we express the higher Johnson homomorphisms as coboundary maps in group cohomology.
Our methods involve the topology of nilpotent homogeneous spaces and lattices in nilpotent Lie algebras.
In particular, we develop a notion of the ``polynomial straightening'' of a singular homology chain in a nilpotent homogeneous space.
\end{abstract}

\section{Introduction}
\subsection{Summary of Results}
Let $\pi$ be a free group of rank $n$.
For each $k>0$, let $\Gamma_k$ be the largest $(k-1)$-step nilpotent quotient of $\pi$.
The \emph{$k$th Andreadakis group} $\Aker{n}{k}$ is the kernel of the action of the automorphism group $\Aut \pi$ on $\Gamma_k$.
The $\{\Aker{n}{k}\}_k$ are a series of nested normal subgroups in $\Aut \pi$.
The \emph{$k$th free group Johnson homomorphism} $\fjh{k}$ is a homomorphism from $\Aker{n}{k}$ to a finite-rank abelian group that encodes the action of $\Aker{n}{k}$ on $\Gamma_{k+1}$.
We give a full definition in Section~\ref{ss:johom} below.

If $n=2g$ is even, we identify $\pi$ with the fundamental group of $\Sigma=\Sigma_{g,1}$, the genus-$g$ surface with a single boundary component.
The \emph{mapping class group} $\bdM=\bdM(\Sigma,\partial \Sigma)$ of $\Sigma$ relative to $\partial \Sigma$ is the group of diffeomorphisms of $\Sigma$ fixing $\partial\Sigma$ pointwise, modulo equivalence by homotopy relative to $\partial\Sigma$.
Since $\Sigma$ is an Eilenberg--Mac~Lane $K(\pi,1)$ space, we have an embedding $\bdM\into \Aut \pi$.
The \emph{$k$th Torelli subgroup} is $\bdI{k}=\Aker{n}{k}\cap\bdM$, the kernel of the action of $\bdM$ on $\Gamma_k$.
The \emph{$k$th surface Johnson homomorphism} $\johom{k}$ is essentially the composition of $\bdI{k}\into \Aker{n}{k}$ with $\fjh{k}$.
The \emph{$k$th Morita homomorphism} $\mohom{k}$ is a related homomorphism from $\bdI{k}$ to $H_3(\Gamma_k)$ defined in terms of the action of $\bdI{k}$ on the group homology chains of $\pi$.
We define these maps in Section~\ref{ss:johom} and Section~\ref{ss:mohom} below.

In this paper we consider extensions of the $\{\delta_k\}_k$ to maps defined on $\Aut \pi$ and extensions of $\{\johom{k}\}_k$ and $\{\mohom{k}\}_k$ to $\bdM$.
Morita initiated this line of inquiry in~\cite{mejo}.
Several authors have contributed to these ideas, constructing extensions of the Johnson homomorphisms with various properties.
These include Hain~\cite[Section~14.6]{hain}, Perron~\cite{perron}, Kawazumi~\cite{kawazumi}, Morita--Penner~\cite{moritapenner}, and Bene--Kawazumi--Penner~\cite{bkp}.
The main motivation is to understand the structure of the \emph{Torelli group} $\mathcal{I}=\bdI{2}$ and the group $\IA_n=\Aker{n}{2}$.
The Torelli group is kernel of the linear representation of $\bdM$ given by the action on $H_1(\Sigma)$ and is the more mysterious part of $\bdM$.
In particular, the study of these extensions could help to identify the images of the $\{\johom{k}\}_k$ and $\{\fjh{k}\}_k$.

In the author's previous paper~\cite{day}, we analyzed maps from $\Sigma\times [0,1]$ into nilpotent homogeneous spaces to construct a crossed homomorphism (a group cohomology $1$-cocycle) on $\bdM$ that extends the $k$th Morita homomorphism for each $k\geq 2$.
The ranges of these maps are finite-dimensional real vector spaces~\cite[Main Theorem~A]{day}.
By composing the extension of the $k$th Morita homomorphism with a particular map, we found a crossed homomorphism extending the $k$th Johnson homomorphism to a map from $\bdM$ to a finite-dimensional vector space~\cite[Main Theorem~B]{day}.
In the present paper, we make several improvements to those results.

Informally, the Johnson homomorphism $\fjh{k}$ is projection of a group cohomology coboundary map with ``nonabelian coefficient module" (see Definition~\ref{de:fjh}).
We express a version of $\fjh{k}$ as a \emph{genuine} coboundary, and in doing so extend it to a crossed homomorphism.
Let $\lag{k}$ be the Lie algebra of the Mal'cev completion of $\Gamma_k$ (see Theorem~\ref{th:malcev} below) and let $\lag{k}^{(1)}$ be its commutator subalgebra.
Let $H=H_1(\pi)$ be the abelianization of $\pi$.
As we explain in the text, there are induced actions of $\Aut \pi$ on $\lag{k}$ and $H$ that extend to objects built functorially out of them.
\begin{theorem}\label{th:mautfncobd}
For each $k\geq 2$, there is an exact sequence
\[0\to \Hom(H,\lag{k}^{(1)})\to \Hom(H,\lag{k}) \to \Hom(H,H\otimes \R)\to 0,\]
of $(\Aut \pi)$--modules that are finite-dimensional real vector spaces, 
and an injective, $(\Aut \pi)$--equivariant homomorphism $\eta\co \Range(\fjh{k})\to \Hom(H,\lag{k}^{(1)})$, 
such that any cocycle representing the cohomology class $d(\mathid\otimes 1)$ is a crossed homomorphism extending $\eta\circ\fjh{k}$.
\end{theorem}
Here $d$ is the coboundary map
\[d\co H^0(\Aut\pi;\Hom(H,H\otimes\R))\to H^1(\Aut\pi; \Hom(H,\lag{k}^{(1)}))\] 
from the long exact sequence for the cohomology of $\Aut \pi$
and $\mathid\otimes1$ is the canonical inclusion $H\into H\otimes\R$.
Since $\mathid\otimes 1$ is invariant under the natural action of $\Aut\pi$ on $\Hom(H,H\otimes\R)$, we consider it as an element of $H^0(\Aut\pi;\Hom(H,H\otimes\R))$.
We prove Theorem~\ref{th:mautfncobd} by using a map of a topological graph into a nilpotent homogeneous space to describe and analyze a lift of $\mathrm{id}\otimes 1$ to $\Hom(H,\lag{k})$.
The proof of Theorem~\ref{th:mautfncobd} appears in Section~\ref{ss:autfnejh} below.
We also provide an interpretation of the surface Johnson homomorphisms as coboundaries in Section~\ref{ss:mcg}.

By using natural lattices in the $\{\Hom(H,\lag{k})\}_k$, we are able to refine the ranges of our maps.
The following corollary is restated as Corollary~\ref{co:nicest} in Section~\ref{ss:range} below.
\begin{corollary}\label{co:frab}
For each $k\geq 2$, there is a finite-rank free abelian group $A_k$ with an $(\Aut \pi)$--action, an injective, $(\Aut \pi)$--equivariant homomorphism
$\eta\co \Range(\fjh{k})\to A_k$,
and a crossed homomorphism $\gamma_k\co \Aut\pi\to A_k$ such that $\gamma_k$ extends $\eta\circ\fjh{k}$.
\end{corollary}
We also deduce a version of this corollary for surface Johnson homomorphisms, Corollary~\ref{co:surnice} below.
In Section~\ref{ss:example}, we explicitly examine the crossed homomorphism from Corollary~\ref{co:frab} in the case $k=2$.
We find that it can be identified with a crossed homomorphism constructed by Morita in~\cite{mejo}.

From Corollary~\ref{co:frab}, we get an immediate structural corollary about $\Aut \pi$.
\begin{corollary}\label{co:semidirect}
For $k\geq 2$, there is an induced action of $(\Aut\pi)/\Aker{n}{k}$ on $A_k$ and there is an injective homomorphism to the semidirect product:
\[
\begin{split}
(\Aut \pi)/\Aker{n}{k+1}&\to ((\Aut\pi)/\Aker{n}{k})\ltimes A_k\\
\alpha\cdot \Aker{n}{k+1} &\mapsto (\alpha\cdot\Aker{n}{k},\gamma_k(\alpha)).\\
\end{split}
\]
\end{corollary}

\begin{proof}
By Lemma~\ref{le:commact} below, $\Aker{n}{k}$ acts trivially on $\Hom(H,\lag{k+1}^{(1)})$ and its submodule $A_k$, so the action $\Aut\pi\curvearrowright A_k$ induces the required action.
It follows from the definitions that the map is a well-defined homomorphism.
Let $\alpha$  be in $\Aut\pi$ with $\alpha\cdot \Aker{n}{k+1}$ in the kernel of the above map.
Then $\alpha$ is in $\Aker{n}{k}$, and is therefore in $\ker (\gamma_k|_{\Aker{n}{k}})=\ker \delta_k=\Aker{n}{k+1}$.
\end{proof}
There is an analogous corollary for $\bdM$ that we leave to the reader to formulate.
These results are an improvement on the similar result in Day~\cite[Main~Corollary~A]{day} since the group $A_k$ is a finite-rank abelian group.

Finally, we build extensions of the Morita homomorphisms $\{\mohom{k}\}_k$, improving on Day~\cite[Main Theorem~A]{day}.
\begin{theorem}\label{th:mobetter}
For each $k\geq 2$, there is a finite-rank free abelian group $B_k$ with a $\bdM$--action, an injective, $(\Aut \pi)$--equivariant homomorphism $\eta\co\Range(\mohom{k})\to B_k$, and a crossed homomorphism $\epsilon_k\co \bdM\to B_k$ such that $\epsilon_k$ extends $\eta\circ\mohom{k}$.
\end{theorem}
To prove this theorem, we develop a theory of polynomial straightenings of simplices in nilpotent homogeneous spaces.
Let $\Delta$ be a simplex in a nilpotent homogeneous space $X$ with all of its vertices mapping to the basepoint $[e]$.
The \emph{polynomial straightening} $s(\Delta)$ of $\Delta$ is a simplex in $X$ in the same homotopy class relative to vertices as $\Delta$, but defined by a polynomial function in exponential coordinates, and such that each edge of $s(\Delta)$ lifts to a left-translation of a $1$-parameter family in the universal cover of $X$ (see Definition~\ref{de:polystr}).

Morita's original, group-theoretic construction of $\mohom{k}$ involves a procedure for building group homology chains that encode data from elements of $\bdM$.
We transform these resulting chains into straightened polynomial chains in a nilpotent homogeneous space,
and then turn these chains into Lie algebra chains using an integration procedure.
The key observation is that the integral of a polynomial with integer coefficients over the standard simplex takes values in a subset of $\Q$ with bounded denominators (see Proposition~\ref{pr:Kexists}).
This appears in Section~\ref{se:moext}.

\subsection{Layout of the paper}
Section~\ref{se:preco} contains background, conventions and preliminary results: the Johnson and Morita homomorphisms are explained in Sections~\ref{ss:johom} and~\ref{ss:mohom}, background on nilpotent Lie groups, Lie algebras and homogeneous spaces is in Sections~\ref{ss:nillie} and~\ref{ss:nilhomog}, and important notation and constructions from Day~\cite{day} appear in Section~\ref{ss:nilstudy}.
Section~\ref{se:extjohom} discusses extensions of Johnson homomorphisms: the proof Theorem~\ref{th:mautfncobd} is in Section~\ref{ss:autfnejh}, an analogous theorem for mapping class groups is in Section~\ref{ss:mcg}, the proof of Corollary~\ref{co:frab} is in Section~\ref{ss:range}, and an example appears in Section~\ref{ss:example}.
Section~\ref{se:moext} is about extensions of Morita homomorphism: Section~\ref{ss:polynomialchains} develops a theory of polynomial homology chains for nilpotent homogeneous spaces, and Section~\ref{ss:mocareful} contains the proof of Theorem~\ref{th:mobetter}.

\subsection{Acknowledgments}
I would like to thank Benson Farb and Danny Calegari for comments on a draft of this paper.
This research was done under the support of an N.S.F. Mathematical Sciences Postdoctoral Research Fellowship.

\section{Preliminaries and conventions}\label{se:preco}
\subsection{Conventions}
All vector spaces,  Lie groups and Lie algebras are over the reals and all tensor products are taken over the reals.
Homology of spaces is singular homology with \emph{piecewise-smooth chains}.
This is so that differential forms can be integrated across singular chains without further comment.
Cohomology of spaces is de Rham cohomology.

\subsection{Johnson homomorphisms}\label{ss:johom}
First we review the terms from the introduction, with precise definitions.
As above, let $\pi$ be a nonabelian free group of rank $n$.
The lower central series $\{\lcs{\pi}{k}\}_k$ of $\pi$ has $\lcs{\pi}{0}=\pi$ and $\lcs{\pi}{k+1}=[\pi,\lcs{\pi}{k}]$.
For each $k\geq 1$, let $\trunc{k}$ be the quotient $\pi/\lcs{\pi}{k-1}$, the \emph{$(k-1)$st nilpotent truncation} of $\pi$.
This is the largest $(k-1)$--step nilpotent quotient of $\pi$.
The lower central series subgroups are characteristic subgroups of $\pi$ so the action $\Aut \pi\curvearrowright\pi$ descends to actions $\Aut\pi\curvearrowright\trunc{k}$.
The kernel of the action of $\Aut\pi$ on $\trunc{k}$ is the \emph{$k$th Andreadakis group} $\Aker{n}{k}$.
The study of the $\{\Aker{n}{k}\}_k$ groups was initiated by Andreadakis~\cite{andreadakis}.

For $\alpha\in\Aker{n}{k}$ and $u\in\pi$, it follows immediately from the definitions that $\alpha(u)u^{-1}$ maps to the trivial element of $\trunc{k}$ and is therefore in $\lcs{\pi}{k-1}$.
Let $H=H_1(\pi)$ be the abelianization of $\pi$ and let $\lv{k+1}$ be the abelian group $\lcs{\pi}{k-1}/\lcs{\pi}{k}$.
\begin{definition}\label{de:fjh}
The $k$th \emph{free-group Johnson homomorphism} $\fjh{k}$ is the homomorphism 
\[
\fjh{k}\co\Aker{n}{k}\to \Hom(H,\lv{k+1})\]
given by the pairing
\[\begin{split}
\Aker{n}{k}\times H&\to \lv{k+1}\\
(\alpha,[u])&\mapsto (\alpha(u)u^{-1})\cdot\lcs{\pi}{k},
\end{split}\] 
for $u\in\pi$.
\end{definition}
As explained in Satoh~\cite[Section~2.4]{satoh}, 
this is a well-defined, $(\Aut \pi)$--equivariant homomorphism.
Informally, $\fjh{k}$ is a projection of a group cohomology coboundary in nonabelian coefficients: the expression $\alpha(u)u^{-1}$ is a nonabelian version of the coboundary expression from equation~\eqref{eq:groupcoboundary} below.

As mentioned above, if $n=2g$, we fix a genus--$g$ surface $\Sigma$ with one boundary component and a basepoint $*$ on the boundary, and identify $\pi$ with $\pi_1(\Sigma,*)$.
Since $\Sigma$ is a $K(\pi,1)$, the action of $\bdM$ on $\pi$ induces an inclusion $\bdM\into\Aut\pi$.
The preimage of $\Aker{2g}{k}$ in $\bdM$ is the \emph{$k$th Torelli group} $\bdI{k}$.
Poincar\'e--Lefschetz duality $D\co H^*\to H$ induces an isomorphism $D_*\co \Hom(H,\lv{k+1})\to H\otimes \lv{k+1}$.
Let $i\co\bdI{k}\to \Aker{2g}{k}$ denote the inclusion.
\begin{definition}\label{de:johom}
The $k$th \emph{Johnson homomorphism} $\johom{k}$ for the bounded surface $\Sigma$ is the map
\[\johom{k}=D_*\circ\fjh{k}\circ i\co \bdI{k}\to H\otimes\lv{k+1}.\]
\end{definition}
This definition is equivalent to the definition appearing in Morita~\cite[Section~2]{mejo}.
It is customary to think of the range of the Johnson homomorphism for a surface as being $H\otimes\lv{k+1}$ instead of $\Hom(H,\lv{k+1})$.
There are also theoretical advantages to this point of view, as seen in Theorem~\ref{th:morita's} below.
Johnson first defined Johnson homomorphisms (for surfaces) in~\cite{joabquo}; higher Johnson homomorphisms appeared in Johnson~\cite{josurv} and in Morita~\cite{mejo}.
The study of Johnson homomorphisms for free groups came later, starting with work of Kawazumi~\cite{kawazumi}.

\subsection{Homology of groups}
Let $\Delta_m$ be the standard $m$--simplex, with vertices labeled $v_0,\dotsc, v_m$.
Let $G$ be a group.
A \emph{$G$--labeling} of $\Delta_m$ is a labeling of each directed edge of $\Delta_m$ by an element of $G$, such that the product of the labels, reading around the boundary of any triangle face of $\Delta_m$ in the same direction, is the trivial element.
A unique labeling of the $m$--simplex $\Delta_m$ is induced by a labeling of only the $m$ edges on the path $v_0-v_1-\dotsb-v_m$.
The standard group homology chains of $G$ are the complex with $C_m(G)$ the free abelian group generated by the set of all $G$--labeled $m$--simplices.
The complex $C_*(G)$ has the usual boundary operation for simplices, while remembering labels.
Reading off the labels on the path $v_0-v_1-\dotsb-v_m$ induces an identification of the basis simplices of $C_m(G)$ with the $m$--fold Cartesian product of $G$.
This identifies $C_*(G)$ as defined here with the \emph{standard complex} for group homology, as appears in Chapter~I.5 of Brown~\cite{brown}.
The \emph{group homology} of $G$ is $H_*(G)=H_*(C_*(G))$.
Note that a map $G_1\to G_2$ induces maps $C_*(G_1)\to C_*(G_2)$ and $H_*(G_1)\to H_*(G_2)$ by acting on labels.

Let $G$ be a group and $A$ a $G$--module.
A \emph{crossed homomorphism} from $G$ to $A$ is a group cohomology $1$--cocycle in $Z^1(G;A)$, that is, a map $f\co G\to A$ satisfying
\[f(gh)=g\cdot f(h)+f(g), \quad \text{for $g,h\in G$.}\]
A \emph{principal crossed homomorphism} from $G$ to $A$ is a group cohomology $1$--coboundary in $B^1(G;A)$, that is, for fixed $a\in A$, a map $d a\co G\to A$ given by
\begin{equation}\label{eq:groupcoboundary}d a(g)= g\cdot a - a.\end{equation}
The first cohomology group $H^1(G;A)$  is $Z^1(G;A)/B^1(G;A)$.
The zeroth cohomology group $H^0(G;A)=Z^0(G;A)$ is the submodule of invariant elements $A^G$.

\subsection{Morita's homomorphisms}\label{ss:mohom}
In~\cite{moabquo}, Morita defined a series of homomorphisms that refine the Johnson homomorphisms.
We recall that construction here.
Since group homology chains are functorial, the projection $p\co \pi\to\trunc{k}$ induces a map $p_*\co C_*(\pi)\to C_*(\trunc{k})$, and the action $\bdM\curvearrowright\pi$ induces an action $\bdM\pi\curvearrowright C_*(\pi)$.
With the induced orientation, the boundary loop $\partial\Sigma$ defines an element $\ell\in\pi$.
Let $[\ell]\in C_1(\pi)$ denote the $1$--simplex labelled by $\ell$.
The action of $\bdM$ fixes $\ell$ and therefore also fixes $[\ell]$.
Let $C_\Sigma\in C_2(\pi)$ have $\partial C_\Sigma=[\ell]$, with $C_\Sigma$ mapping to a representative of the generator of $H_2(\pi/\overline{\langle \ell\rangle})$, where $\overline{\langle \ell\rangle}$ is the subgroup normally generated by $\ell$ (basically, $C_\Sigma$ is a fundamental class for $\pi$ relative to $\ell$).
For $[\phi]\in\bdM$, let $D_{[\phi]}\in C_3(\Sigma)$ with $\partial D_{[\phi]}=[\phi]_*C_\Sigma-C_\Sigma$ ($D_{[\phi]}$ exists because $\partial ([\phi]_*C_\Sigma-C_\Sigma)=\ell-[\phi]_*\ell=0$ and $Z_2(\pi)=B_2(\pi)$).
Note that $\partial p_*D_{[\phi]}=[\phi]_*p_*C_\Sigma-p_*C_\Sigma$; so if $[\phi]\in\bdI{k}$ then $\partial p_*D_{[\phi]}=0$.

\begin{definition}\label{de:mohom}
The $k$th Morita homomorphism $\mohom{k}\co \bdI{k}\to H_3(\trunc{k})$ is the map $[\phi]\mapsto [p_*D_{[\phi]}]$, where $\partial D_{[\phi]}=[\phi]_*C_\Sigma-C_\Sigma\in C_2(\pi)$ and $\partial C_\Sigma=[\ell]$ as above.
\end{definition}
The original definition uses $-[\ell]$ instead of $[\ell]$ and $C_\Sigma-\phi_*C_\Sigma$ instead of $\phi_*C_\Sigma-C_\Sigma$, so this definition is equivalent.

This map is connected to the Johnson homomorphism by Morita's theorem~\cite[Theorem~3.1]{moabquo}, restated below.
From the definitions, we have a central extension:
\begin{equation}\label{eq:gpseq}
1\to \lv{k+1}\to \trunc{k+1}\to\trunc{k}\to 1.
\end{equation}
Let $d^2\co H^3(\trunc{k})\to H^1(H;H^1(\lv{k+1}))\cong H\otimes \lv{k+1}$ be the differential from the Hochschild--Serre spectral sequence for the integral cohomology of sequence~\eqref{eq:gpseq}.
From~\cite[Theorem~3.1]{moabquo}:
\begin{theorem}[Morita]\label{th:morita's}
We have $\johom{k}=d^2\circ\mohom{k}$.
\end{theorem}

\subsection{Properties of nilpotent Lie groups}\label{ss:nillie}
The following two theorems are well known; a reference is Raghunathan~\cite[Chapter~2]{raghunathan}.
A \emph{lattice} in a Lie group or vector space is a discrete subgroup with finite covolume.
\begin{theorem}[Mal'cev]\label{th:malcev}
Every finitely generated, torsion-free nilpotent group $\Gamma$ embeds as a lattice in a unique simply connected, nilpotent Lie group $G$, the \emph{Mal'cev completion} of $\Gamma$.
\end{theorem}

\begin{theorem}
\label{th:malcevfunc}
Every map between finitely generated, torsion-free nilpotent groups extends uniquely to a map between their Mal'cev completions.
In particular, if $\Gamma$ is a finitely generated, torsion-free nilpotent group and $G$ is its completion, then $\Aut \Gamma$ embeds in $\Aut G$.
\end{theorem}

Let $G$ be a simply connected, nilpotent Lie group with Lie algebra $\lag{}$.
Recall the Lie group exponential map $\exp\co \lag{}\to G$ (which sends $\vct x\in\lag{}$ to $f(1)$, where $f\co \R\to G$ is a homomorphism with $\frac{\partial}{\partial t}f|_0=\vct x$).
From Corwin--Greenleaf~\cite[Theorem~1.2.1]{corwingreenleaf}:
\begin{theorem}
The map $\exp\co\lag{}\to G$ is an analytic diffeomorphism.
\end{theorem}
This means that we can use the logarithm map $\log=\exp^{-1}\co G\to \lag{}$ as a coordinate system for $G$.
This coordinate system is called \emph{exponential coordinates}.

The \emph{structure constants} of a Lie algebra $\lag{}$ with respect to the basis $\{\vct{x}_1,\dotsc,\vct{x}_n\}$ are the real numbers $c_{i,j}^k$ such that
\[[\vct{x}_i,\vct{x}_j]=\sum_{l=1}^kc_{i,j}^k\vct{x}_k.\]
The next result is from Raghunathan~\cite[Theorem 2.12]{raghunathan}, and the remarks following that theorem.
\begin{theorem}\label{th:lattice}
If $\Gamma$ is a lattice in $G$, then the $\Z$--span of $\log(\Gamma)$ is a lattice $\latt$ in the vector space $\lag{}$.
Further, the structure constants of $\lag{}$ with respect to any basis for $\latt$ are rational.
\end{theorem}

The following is from Corwin--Greenleaf~\cite[Theorem~1.2.1]{corwingreenleaf} and the discussion preceding that theorem.
Let $*\co \lag{}\times\lag{}\to\lag{}$ denote the pullback of the Lie group product to $\lag{}$:
\[\vct{x}*\vct{y}=\log(\exp(\vct{x})\exp(\vct{y})),\quad \text{for $\vct{x},\vct{y}\in\lag{}$}.\]
\begin{theorem}
There is a universal formula for $\vct{x}*\vct{y}$ as a rational linear combination of nested commutators of $\vct{x}$ and $\vct{y}$, that is satisfied for any pair of elements $\vct{x},\vct{y}$ in any nilpotent Lie algebra $\lag{}$.

\end{theorem}
This formula is called the \emph{Baker--Campbell--Hausdorff formula}.
The Baker--Campbell--Hausdorff formula is an infinite sum of rational multiples of nested commutators of the terms $\vct{x}$ and $\vct{y}$; since we are working in a nilpotent Lie algebra, all but finitely many of these terms are zero and the formula is a finite sum.
Our discussion would gain little from the full expression for this formula, but the interested reader can find it in Corwin--Greenleaf~\cite[p.\,11]{corwingreenleaf}.
Instead we display the first few terms:
\begin{equation}\vct{x}*\vct{y}= \vct{x}+\vct{y}+\frac{1}{2}[\vct{x},\vct{y}]+\frac{1}{12}[\vct{x},[\vct{x},\vct{y}]]+\frac{1}{12}[\vct{y},[\vct{x},\vct{y}]]+\dotsb,\end{equation}
where the ellipses conceal the commutators of four or more terms.

By a polynomial map between vector spaces, we mean a map whose coordinate functions are polynomials in the coordinates of the domain.
This notion, and the polynomial degree of such a map, are independent of any choice of coordinates.
\begin{corollary}\label{co:polymult}
The map $*\co\lag{}\times\lag{}\to \lag{}$ is a polynomial map whose degree is bounded by the nilpotence class of $\lag{}$.
\end{corollary}

\begin{proof}
Suppose $\{\vct{x}_1,\dotsc,\vct{x}_n\}$ is a basis for $\lag{}$.
If $\vct{x}=\sum r_i \vct{x}_i$ and $\vct{y}=\sum s_i \vct{x}_i$, then any nested commutator in $\vct{x}$ and $\vct{y}$ expands into a linear combination of nested commutators of the $\{\vct{x}_i\}_i$, with coefficients being monomials in $\{r_i,s_i\}_i$.
Of course, nested commutators in $\{\vct{x}_i\}_i$ evaluate into linear combinations of $\{\vct{x}_i\}_i$.
The degree of a commutator's coefficient as a monomial in $\{r_i,s_i\}$ is the depth of the nested commutator.
If a given nested commutator evaluates to a nonzero vector, then the degree of its coefficient is bounded by the nilpotence class of $\lag{}$.
Then by expanding and evaluating the commutators in the Baker--Campbell--Hausdorff formula, we expresses $\vct{x}*\vct{y}$ as a linear combination of the $\{\vct{x}_i\}_i$, where the coefficients are polynomials in the $\{r_i,s_i\}_i$, with the stated bound on degree.
\end{proof}

In fact, we can say something a little stronger about the multiplication map as a polynomial.
We need the following fact from Corwin--Greenleaf~\cite[Theorem~5.4.2]{corwingreenleaf}.
\begin{theorem}\label{th:latticelattice}
Let $\Gamma$ be a lattice in $G$.
Then $\Gamma$ is contained as a finite-index subgroup in a lattice $\Gamma_1$ in $G$, such that the image $\log(\Gamma_1)$ is a lattice in $\lag{}$.
\end{theorem}

\begin{corollary}\label{co:multlattice}
If $\Gamma$ is a lattice in $G$, then there is a lattice $\latt$ in $\lag{}$ containing $\log(\Gamma)$ such that $*\co \lag{}\times\lag{}\to \lag{}$ maps the lattice $\latt\times\latt$ to $\latt$.
\end{corollary}

\begin{proof}
Take $\latt$ to be $\log(\Gamma_1)$ as in Theorem~\ref{th:latticelattice}.
Then each element of $\latt\times\latt$ is of the form $(\log(g_1),\log(g_2))$ for $g_1,g_2\in \Gamma_1$.
Since $\Gamma_1$ is a subgroup, $*$ sends this element to $\log(g_1g_2)\in\latt$.
\end{proof}

\subsection{Nilpotent Lie algebras and homogeneous spaces}\label{ss:nilhomog}
Let $\Gamma$ be a torsion-free, finitely generate nilpotent group and $G$ its Mal'cev completion.
Let $\lag{}$ be the Lie algebra of $G$.
The \emph{standard Lie algebra chain complex} $C_*(\lag{})$ has chain groups $C_k(\lag{})=\bigwedge^k\lag{}$ (the exterior power) and a boundary map determined by the Lie bracket:
\begin{equation}\label{eq:lieboundary}
\begin{split}
\partial(\vct{x_1}&\wedge\dotsb\wedge\vct{x_k})=\\
&\sum_{1\leq i<j\leq k}(-1)^{i+j+1}[\vct{x}_i,\vct{x}_j]\vct{x}_1\wedge\dotsb\wedge\hat{\vct{x}}_i\wedge\dotsb\wedge\hat{\vct{x}}_j\wedge \dotsb\wedge \vct{x}_k,
\end{split}
\end{equation}
where $\hat{\vct{x}}_i$ and $\hat{\vct{x}}_j$ indicate the omission of these terms.
The cochain complex $C^k(\lag{})$ has chain groups $C^k(\lag{})=\bigwedge^k\lag{}^*$, the exterior power of the dual space, with a coboundary dual to the boundary described above.
A reference is Weibel~\cite[Chapter 7]{weibel}.

Let $X=G/\Gamma$ be the homogeneous space.
Since $G$ is the Mal'cev completion of $\Gamma$, the group $\Aut \Gamma$ acts on $G$.
The following observation is important, but the proof is straightforward and is omitted.
\begin{proposition}\label{pr:obs}
The action $\Aut\Gamma\curvearrowright G$ induces a smooth action on $X$, and
the induced action on $\pi_1(X,[e])\cong \Gamma$ is the usual one.
\end{proposition}
This action gives us a notion of left-invariant differential form on $X$.
In fact, each member of $\bigwedge^k \lag{}^*$ defines a left invariant form on $X$: the unique left-invariant form that restricts to that functional on the tangent space to $[e]$.
Let $\lprop\co C^k(\lag{})\to C^k(X;\R)$ denote this map, which is a chain homomorphism.
Let $\langle,\rangle$ denote the evaluation pairing $C^*(\lag{})\otimes C_*(\lag{})\to\R$.
From Day~\cite[Definition~2.2]{day}, we have an adjoint version of $\lprop$:
\begin{definition}\label{de:vectorize}
Let $v\co C_*(X;\R)\to C_*(\lag{})$ send 
$C\in C_n(X;\R)$ to the unique multi-vector $v(C)\in C_n(\lag{})$ such that 
\[\int_C \lprop(\alpha)=\langle\alpha,v(C)\rangle,\]
for every $\alpha\in C^n(\lag{k}).$
\end{definition}

A theorem of Nomizu~\cite{nomizu} states that $\lprop$ induces an isomorphism from the Lie algebra cohomology of $\lag{}$ to the de Rham cohomology of $X$.
In~\cite[Corollary~2.1]{day} we proved the following as a corollary to Nomizu's theorem. 
\begin{corollary}\label{co:vecisom}
The induced map $v_*\co H_*(X;\R)\to H_*(\lag{})$ is an isomorphism.
\end{corollary}

We state a few more properties of $v$.
The next one follows immediately from the definition and needs no proof.
\begin{proposition}\label{pr:basisform}
Let $B$ be a basis for $C_n(\lag{k})$ and let $B^*$ be its evaluation-dual basis for $C^n(\lag{k})$.
For $\vct{x}\in B$, let $\vct{x}^*\in B^*$ be its dual.
Then for any piecewise smooth $C\in C_n(X)$,
\[v(C)=\sum_{\vct{x}\in B}\biggl(\int_C \lprop(\vct{x}^*)\biggr)\vct{x}.\]
\end{proposition}

\begin{lemma}\label{le:dim}
Let $Y$ be an $m$--dimensional simplicial complex with smooth structures on its $m$--simplices, and let $f\co Y\to X$ be a piecewise-smooth map.
If $l>m$ and $C\in C_l(Y)$, then $v(f_*C)=0$.
\end{lemma}

\begin{proof}
By Proposition~\ref{pr:basisform}, $v(f_*C)$ is a linear combination of elements in a basis $B$ for $C_l(\lag{})$ with coefficients given as integrals of $l$--forms across $f_*C$.
For $\alpha\in B^*$, we use the pullback to compute $\int_{f_*C}\alpha=\int_Cf^*\alpha$.
However, since $Y$ is $m$--dimensional, all $l$--forms are degenerate on $Y$ and $\int_Cf^*\alpha=0$.
\end{proof}

\begin{lemma}\label{le:fc}
Suppose $T_0$ and $T_1$ are two different fundamental classes relative to the boundary for an $m$--manifold $M$ possibly with boundary, and suppose $f\co M\to X$ is piecewise smooth.
Then $v(f_*T_1)=v(f_*T_0)\in C_m(\lag{})$.
\end{lemma}

\begin{proof}
The hypotheses imply that there is are chains $C\in C_{m+1}(M)$ and $C'\in C_m(\partial M)$ with $\partial C=C'+T_1-T_0$.
Lemma~\ref{le:dim} implies that $v(f_*C)=0$ and $v(f_*C')=0$ (since $\partial M$ is $(m-1)$--dimensional).
Then $0=\partial v(f_*C)= v(f_*T_1)-v(f_*T_0)$.
\end{proof}

By Theorem~\ref{th:malcevfunc}, each homomorphism $\Z\to\Gamma$ extends uniquely to a homomorphism $\R\to G$.
Therefore we can define the following.
\begin{definition}\label{de:canonical}
For $\gamma\in \pi_1(X,[e])=\Gamma$, let the \emph{canonical representative} of $\gamma$ be the unique map $\R/\Z\to G/\Gamma$ that lifts to a homomorphism $\R\to G$ sending $1$ to $\gamma$.
\end{definition}

\begin{lemma}\label{le:vecabel}
If $G$ is abelian, then the maps
\[\Gamma\into G \xrightarrow{\log} \lag{} \quad \text{and}\]
\[\Gamma \cong H_1(X) \xrightarrow{v_*} H_1(\lag{})\cong\lag{}\]
are identical.
\end{lemma}

\begin{proof}
Let $c\co [0,1]\to X$ be the canonical representative of $\gamma\in \pi_1(X,[e])$, and let $\tilde c\co [0,1]\to G$ be its lift to $G$ starting at $e$.
Then $\log(\gamma)$ is the (one-sided) derivative of $\tilde c$ at $0$.
Of course, $c$ represents $[\gamma]$ in $H_1(X)$.
In exponential coordinates, $\tilde c$ is a linear map, and since $G$ is abelian, the left-invariant differential forms are all constant.
This means that the integral of any left-invariant differential $1$--form across $\tilde c$ will equal its value on the derivative of $\tilde c$ at $e$.
By Definition~\ref{de:vectorize}, this proves the lemma.
\end{proof}

\subsection{Polynomials respecting lattices}
This material is important for our refinements of the ranges of crossed homomorphisms.
\begin{proposition}\label{pr:integervaluedcoeffs}
Fix positive integers $n$ and $d$.
There is a positive integer $K$ depending on $n$ and $d$, such that
for any polynomial map $f\co \R^n\to \R$ of degree $d$ sending $\Z^n$ into $\Z$,
the coefficients of $f$ are in $\frac{1}{K}\Z$.
\end{proposition}

\begin{proof}
This follows from Prasolov~\cite[Theorem~3.2.3]{prasolov}, which states that any polynomial $f$ from $\R^n$ to $\R$, sending $\Z^n$ to $\Z$, is an integer linear combination of products of binomial coefficient polynomials $x_i\mapsto \binom{x_i}{k}$ for various $k\in\Z$ ($x_i$ is a coordinate variable for $\R^n$).
\end{proof}

On a vector space $V$, we identify the tangent space $T_0V$ with $V$ by the exponential map.
The derivatives of the translation functions canonically identify the tangent spaces for different points.
Therefore there is a canonical identification of the tangent bundle $TV$ to $V\times V$ (the first coordinate denoting the point and the second the vector).
We tacitly use this identification in later arguments.

\begin{lemma}\label{le:polyderiv}
Suppose $f\co V_1\to V_2$ is a polynomial map of vector spaces, and suppose $A_i<V_i$ are lattices such that $f(A_1)\subset A_2$.
There is an integer $K$, depending only on $\dim V_1$ and the degree of $f$, such that the total derivative $Df\co TV_1\to TV_2$ is a polynomial map from $V_1\times V_1\to V_2\times V_2$ that sends $A_1\times A_1$ into $A_2\times \frac{1}{K}A_2$.
\end{lemma}

\begin{proof}
Recall that the total derivative of a map is the induced map on tangent bundles.
The tangent bundle of a vector space is naturally diffeomorphic to its Cartesian product with itself.
The map $Df\co V_1\times V_1\to V_2\times V_2$ sends $(\vct{x},\vct{y})$ to $(f(\vct{x}),\sum y_i(\frac{\partial}{\partial x_i}f)|_{\vct{x}})$.
By picking bases for $A_1$ and $A_2$, we get coordinate polynomials for $f$ that are integer-valued polynomials.
By Proposition~\ref{pr:integervaluedcoeffs}, there is a $K$ such that the coefficients of the coordinate polynomials of $f$ are in $\frac{1}{K}\Z$.
Of course, all partial derivatives of the coordinate polynomials of $f$ also have coefficients in $\frac{1}{K}\Z$.
Then the second-coordinate map of $Df$ sends $A_1$ into $\frac{1}{K}A_2$.
\end{proof}

\begin{lemma}\label{le:polycomp}
Suppose $f\co V_1\to V_2$ and $g\co V_2\to V_3$ are polynomial maps between vector spaces.
Suppose $A_i<V_i$ are lattices and suppose there is an integer $K_1>0$ such that $f(A_1)\subset \frac{1}{K_1}A_2$ and $g(A_2)\subset A_3$.
Then there is an integer $K_2>0$ depending only on $K_1$, $\dim V_1$, $\dim V_2$, $\deg f$ and $\deg g$, such that $g\circ f(A_1)\subset \frac{1}{K_2}A_3$.
\end{lemma}

\begin{proof}
Pick bases for each of the $A_i$.
This gives us coordinate polynomials for $f$ and $g$.
By Proposition~\ref{pr:integervaluedcoeffs}, there are integers $K_3,K_4>0$, depending only on $K_1$, $\dim V_1$, $\dim V_2$, $\deg f$ and $\deg g$, such that the coefficients of the coordinate polynomials for $f$ are in $\frac{1}{K_3}\Z$ and the coefficients for the coordinate polynomials of $g$ are in $\frac{1}{K_4}\Z$.
Set $K_2=K_4\cdot K_3^{\deg g}$.
Then the coefficients for the coordinate functions of the compositions $g\circ f$ are in $\frac{1}{K_2}\Z$, and $g\circ f$ sends $A_1$ into $\frac{1}{K_2}A_3$.
\end{proof}

\subsection{Constructions related to Johnson and Morita homomorphisms} \label{ss:nilstudy}
In Day~\cite{day}, we used nilpotent homogeneous spaces and nilpotent Lie algebras to study Johnson and Morita homomorphisms.
In this section we recall some of those constructions and prove some general lemmas.

Since each $\trunc{k}$ is a torsion-free, finitely-generated nilpotent group, it embeds as a lattice in its Mal'cev completion $G_k$.
Let $X_k$ be the homogeneous space $G_k/\trunc{k}$; this is a compact manifold (an iterated torus bundle) and a $K(\trunc{k},1)$ space.
Let $C_*(X_k;\R)$ denote the piecewise-smooth simplicial chains on $X_k$.
Let $\lag{k}$ denote the Lie algebra of $G_k$.
We use the maps $L\co C^*(\lag{k})\to C^*(X_k;\R)$ and $v\co C_*(X_k)\to C_*(\lag{k})$ as defined above.

Sequence~\eqref{eq:gpseq} induces a central extension of Mal'cev completions:
\[1\to \lv{k+1}\otimes\R \to G_{k+1}\to G_k\to 1.\]
Let $\lal{k+1}$ denote the Lie algebra of $\lv{k+1}\otimes\R$ and let $T_{k+1}$ be the quotient space \mbox{$\lv{k+1}\otimes\R$}$/\lv{k+1}$.
This $T_{k+1}$ is a torus.
Then we have a central extension of Lie algebras
\begin{equation}\label{eq:Liealgs}1\to \lal{k+1}\to \lag{k+1}\to\lag{k}\to 1\end{equation}
and a fiber bundle of homogeneous spaces (with trivial holonomy):
\begin{equation}\label{eq:fibration}T_{k+1}\to X_{k+1}\to X_k.\end{equation}
Since $\Aut \pi$ acts on $\trunc{k}$, it acts on $G_k$, $\lag{k}$, $X_k$, and any objects functorially determined by these.
The maps in the above sequences are all $(\Aut\pi)$--equivariant.
Further, since $\Aker{n}{k}$ acts trivially on $\trunc{k}$, it acts trivially on $X_k$ and~$\lag{k}$.
Of course, $\bdM$ acts on these objects via the inclusion $\bdM\into\Aut\pi$.

\begin{lemma}\label{le:commact}
The group $\Aker{n}{k}$ acts trivially on $\lag{k+1}^{(1)}$ (the commutator subalgebra).
\end{lemma}

\begin{proof}
Let $\alpha\in\Aker{n}{k}$ and $\vct{x},\vct{y}\in\lag{k+1}$; then $\alpha\cdot \vct{x}- \vct{x}$ and $\alpha\cdot \vct{y}-\vct{y}$ map to the trivial element of $\lag{k}$ and are therefore in $\lal{k+1}$. 
Then \[\alpha\cdot [\vct{x},\vct{y}]=[\vct{x}+(\alpha\cdot \vct{x}-\vct{x}),\vct{y}+(\alpha\cdot \vct{y}-\vct{y})]=[\vct{x},\vct{y}]\]
since $\lal{k+1}$ is central in $\lag{k+1}$.
\end{proof}

Now we recall the main construction from Day~\cite{day}.
Let $\rho\co \bdM\curvearrowright X_k$ denote the action from Proposition~\ref{pr:obs}.
It follows from Theorem~\ref{th:malcevfunc} that there is a canonical loop in $X_k$ representing the image in $\pi_1(X_k,[e])=\trunc{k}$ of the boundary class $[\partial\Sigma]\in\pi_1(\Sigma,*)=\pi$.
Fix a piecewise-smooth map $j\co (\Sigma,*)\to (X_k,[e])$ inducing the canonical projection $\pi\to \trunc{k}$ on fundamental groups and sending $\partial\Sigma$ to the canonical representative of the image of $[\partial\Sigma]$.
For $\phi\in\Diff(\Sigma,\partial\Sigma)$, we consider the map $\rho([\phi])\circ j\circ \phi^{-1}$, which is another piecewise-smooth map inducing the canonical projection.
This means there is a piecewise-smooth homotopy $F\co \Sigma\times [0,1]\to X_k$ from $j$ to $\rho([\phi])\circ j\circ\phi^{-1}$ relative to $\partial\Sigma$.
\begin{definition}\label{de:emomap}
For $k\geq 2$, the $k$th extended Morita map $\epsilon_k\co \bdM\to C_3(\lag{k})/B_3(\lag{k})$ is the map
\[[\phi]\mapsto v(F_*[\Sigma\times[0,1]])+B_3(\lag{k}).\]
\end{definition}
Here $[\Sigma\times[0,1]]$ denotes any fundamental class relative to boundary.
We showed in Day~\cite{day} that this is a well-defined crossed homomorphism extending the composition of $\mohom{k}$ with an injective, $\bdM$--equivariant homomorphism.

\section{Extending the Johnson homomorphisms}\label{se:extjohom}
\subsection{Results for $\Aut \pi$}\label{ss:autfnejh}
Let $R$ be the topological graph with $1$ vertex $*$ and $n$ edges, and identify $\pi_1(R,*)$ with $\pi$.
We give $R$ a smooth structure away from $*$.
Choose a piecewise-smooth map $i\co R\to X_{k+1}$ that induces the canonical projection $\pi_1(R,*)\to \pi_1(X_{k+1},[e])$ (the projection $\pi\to\trunc{k+1}$).
\begin{definition}
Let $\jhob=\jhob(i)$ be the function 
\[[c]\mapsto v(i_*c)\]
in $\Hom(H,\lag{k+1})$, where $c\in Z_1(R)$ is any piecewise-smooth representative of $[c]\in H_1(R)$.
\end{definition}

\begin{lemma}
The function $\jhob$ is well defined.
\end{lemma}

\begin{proof}
Suppose $c, c'\in Z_1(R)$ are piecewise-smooth representatives of $[c]\in H_1(R)$.
There is a piecewise-smooth $C\in C_2(R)$ with $\partial C=c-c'$.
By Lemma~\ref{le:dim}, we know that $v(i_*C)=0$.
Since $v$ is a chain map, we then have that $v(i_*c)=v(i_*c')$.
It is easy to see that $\jhob$ is a homomorphism.
\end{proof}

\begin{lemma}\label{le:bdtriv}
The boundary $\partial\co C_2(\lag{k+1})\to B_1(\lag{k+1})$ restricts to the zero map on the kernel of the map $C_2(\lag{k+1})\to C_2(\lag{k})$ induced by the projection in sequence~\eqref{eq:Liealgs}.
\end{lemma}

\begin{proof}
In the standard Lie algebra chain complex, $C_2(\lag{k+1})= \bigwedge^2\lag{k+1}$ and  $C_1(\lag{k+1})=\lag{k+1}$.
The boundary $\partial\co\bigwedge^2\lag{k+1}\to\lag{k+1}$ is given by the commutator map $\vct{x}\wedge \vct{y}\mapsto [\vct{x},\vct{y}]$.
The algebra $\lag{k+1}$ has a basis $B$ such that $B\cap\lal{k+1}$ is a basis for $\lal{k+1}$.
Then $\bigwedge^2\lag{k+1}$ has a basis whose elements are the pairwise wedge products of these basis elements.
It is immediate that the kernel of $C_2(\lag{k+1})\to C_2(\lag{k})$ is generated by the elements of the form $\vct{x}\wedge \vct{y}$ with $\vct{x},\vct{y}\in B$ and $\vct{y}\in\lal{k+1}$.
Since $\lal{k+1}$ is central in $\lag{k+1}$, the generators of the kernel all have trivial boundary.
\end{proof}

The group $\Aut \pi$ acts on $\Hom(H,\lag{k+1})$ as follows: for $f\in\Hom(H,\lag{k+1})$ and $\alpha\in\Aut \pi$, we have $(\alpha\cdot f)(c)=\alpha\cdot (f(\alpha^{-1}\cdot c))$ for $c\in H$.
So we can consider the principal crossed homomorphism $d\jhob\in Z^1(\Aut \pi;\Hom(H,\lag{k+1}))$:
\[d\jhob\co \alpha\mapsto (\alpha\cdot\jhob -\jhob).\]

From Definition~\ref{de:vectorize}, we have a map $v\co C_*(T_{k+1})\to C_*(\lal{k+1})$.
By Corollary~\ref{co:vecisom} this gives us an isomorphism $v_*\co H_*(T_{k+1};\R)\to H_*(\lal{k+1})$.
Let $\vind\co \Hom(H,\lv{k+1})\to \Hom(H,\lag{k+1})$ be given by post-composition with 
\[\lv{k+1}\into \lv{k+1}\otimes\R\cong H_1(T_{k+1};\R)\xrightarrow{v_*}H_1(\lal{k+1})\cong \lal{k+1}\into \lag{k+1}.\]
Note that $\vind$ is an injective, $(\Aut\pi)$--equivariant homomorphism.

Now we have everything needed to cast the Johnson homomorphism as a coboundary.
\begin{keylemma}\label{kl:autfnmain}
For $\alpha\in\Aker{n}{k}$, we have
\[d\jhob(\alpha)=\vind(\fjh{k}(\alpha)).\]
\end{keylemma}

\begin{proof}
Let $\rho\co\Aut \pi\curvearrowright X_{k+1}$ denote the action (from Proposition~\ref{pr:obs}).
Let $T$ be a fixed piecewise-smooth representative of the fundamental class of $S^1$.
Let $\gamma\co (S^1,*)\to (R,*)$ be piecewise-smooth, and let $\bar \gamma$ be $\gamma$ in reverse.
Let $\alpha\in \Aker{n}{k}$.
Since $\alpha\in \Aker{n}{2}$, we know that $\alpha^{-1}\cdot \gamma_*[T]=\gamma_*[T]\in H$.
Then $d\jhob(\alpha)(\gamma_*[T])=\alpha\cdot v(i_*\gamma_*T)-v(i_*\gamma_*T)$.
This, in turn, equals $v(((\rho(\alpha)\circ i\circ\gamma)\cdot (i\circ \bar\gamma))_*T)$, where the $\cdot$ indicates concatenation.

Let $p\co X_{k+1}\to X_k$ denote the projection from the fibration~\eqref{eq:fibration}.
The action of $\Aut \pi$ commutes with the projection; since $\Aker{n}{k}$ acts trivially on $X_k$, we know that $p\circ i\circ \gamma$ and $p\circ \rho(\alpha)\circ i\circ\gamma$ are the same map.

Let $E\subset X_{k+1}$ denote the set of points lying above $p(i(R))$.
Then we have a fiber bundle $T_{k+1}\to E\to p(i(R))$.
Since $p\circ i$ is piecewise smooth, $p(i(R))$ is a finite topological graph.
Since $p\circ i\circ\gamma=p\circ \rho(\alpha)\circ i\circ\gamma\co S^1\to p(i(R))$, we know that $[(\rho(\alpha)\circ i\circ\gamma)\cdot (i\circ \gamma)^{-1}]$ is in the kernel of $\pi_1(E,[e])\to\pi_1(p(i(R)),*)$.
Then there is a homotopy in $E$ relative to $*$ from $i\circ\gamma\cdot (\rho(\alpha\circ i\circ\gamma))^{-1}$ to a loop $\beta\co (S^1,*)\to (T_{k+1},[e])$.
This implies there is a piecewise smooth chain $C\in C_2(E)$ with $\partial C=\rho(\alpha)_*i_*\gamma_*T-i_*\gamma_*T-\beta_*T$.

Since $p_*C$ is supported on $p(i(R))$, Lemma~\ref{le:dim} implies that $v(p_*C)=0$.
Then the naturality of $v$ implies that $v(C)$ is in the kernel of $C_2(\lag{k+1})\to C_2(\lag{k})$.
By Lemma~\ref{le:bdtriv}, $\partial v(C)=0$.
Then since $v$ is a chain map, this means that 
\[0=\partial v(C)= v(\rho(\alpha)_*i_*\gamma_*T-i_*\gamma_*T)-v(\beta_*T).\]
So $v(\beta_*T)=d\jhob(\alpha)(\gamma_*[T])$.
But $[\beta]$ is $\alpha([\gamma])[\gamma]^{-1}\in \pi_1(T_{k+1})=\lv{k+1}$.
So $[\beta]\in\lv{k+1}$ is $\fjh{k}(\alpha)(\gamma_*[T])$.
Since $v(\beta_*T)\in\lag{k+1}$ is the same as $v_*(\beta_*[T])\in H_1(\lal{k+1})=\lal{k+1}\subset\lag{k+1}$, this proves the lemma.
\end{proof}

As usual, $\lag{k+1}^{(1)}$ denotes the commutator subalgebra of $\lag{k+1}$.
The abelianization of $\lag{k+1}$ is $\lag{2}$.
The map $v_*$ induces the following injection:
\begin{equation}H\cong H_1(X_2) \into H_1(X_2;\R) \xrightarrow{v_*} H_1(\lag{2})\cong \lag{2}.\label{eq:Hinj} \end{equation}
This induces an isomorphism $\lag{2}\cong H\otimes\R$, and we have the following exact sequence:
\[0\to \lag{k+1}^{(1)}\to\lag{k+1}\to H\otimes\R \to 0.\]
Since these are free abelian groups, we then have:
\[0\to \Hom(H,\lag{k+1}^{(1)})\to\Hom(H,\lag{k+1})\to \Hom(H,H\otimes\R)\to 0.\]
Then the long exact sequence for the cohomology of $\Aut \pi$ with coefficients in this sequence has the following map:
\[d\co H^0(\Aut \pi;\Hom(H,H\otimes \R))\to H^1(\Aut \pi;\Hom(H,\lag{k+1}^{(1)})).\]

The function $\mathrm{id}\otimes1\in \Hom(H,H\otimes \R)$ sends each element to itself.
Of course, this element is fixed by the action of $\Aut \pi$, so we consider it as $\mathrm{id}\otimes 1\in H^0(\Aut \pi;\Hom(H,H\otimes \R))$.

\begin{proof}[Proof of Theorem~\ref{th:mautfncobd}]
The map $\eta$ from the statement is $\vind$.
It remains to show that any representative of 
\[d(\mathrm{id}\otimes 1)\in H^1(\Aut \pi;\Hom(H,\lag{k+1}^{(1)}))\]
is a crossed homomorphism extending $\vind\circ \fjh{k}
\co \Aker{n}{k}\to\Hom(H;\lal{k+1})$.
Note that $\jhob\in \Hom(H,\lag{k+1})$ is a lift of $\mathrm{id}\otimes 1$.
Key Lemma~\ref{kl:autfnmain} implies that the restriction of $d\jhob$ to $Z^1(\Aut \pi;\Hom(H,\lag{k+1}))$ is a representative of $d(\mathrm{id}\otimes 1)$ that extends $\vind\circ\fjh{k}$.
Any other lift of $\mathrm{id}\otimes 1$ can be written as $\jhob+\eta$, for some $\eta\in\Hom(H,\lag{k+1}^{(1)})$.
By Lemma~\ref{le:commact}, $\Aker{n}{k}$ acts trivially on $\lag{k+1}^{(1)}$ and $d\eta$ restricts to the zero cocycle on $\Aker{n}{k}$.
Then $d(\jhob+\eta)$ extends~$\vind\circ \fjh{k}$.
\end{proof}

We investigate another representative of the class $d(\mathrm{id}\otimes 1)$ in Section~\ref{ss:range} below.

\subsection{Results for mapping class groups}\label{ss:mcg}
In this section we discuss an extension of the surface Johnson homomorphism.
Naturally we assume in this section that $n=2g$ is even.
We start by reviewing a construction from Day~\cite{day}.

The subalgebras $\lag{k+1}^{(1)}$ and $\lal{k+1}$ in $\lag{k+1}$ are both $\bdM$--invariant, so $\lag{k+1}^{(1)}\wedge\lal{k+1}$ is a $\bdM$--invariant subspace of $C^2(\lag{k+1})$.
This gives us a quotient action of $\bdM$ on $C^2(\lag{k+1})/(\lag{k+1}^{(1)}\wedge\lal{k+1})$.
In Day~\cite{day}, we defined a $\bdM$--equivariant homomorphism 
\[\tilde d^2\co C_3(\lag{k})/B_3(\lag{k})\to C_2(\lag{k+1})/(\lag{k+1}^{(1)}\wedge\lal{k+1}).\]
For $C\in C_3(\lag{k})$, this $\tilde d^2([C])$ denotes the class in $C_2(\lag{k+1})/(\lag{k+1}^{(1)}\wedge\lal{k+1})$  of $\partial \tilde C$, where $\tilde C\in C_3(\lag{k+1})$ is any lift of $C$.

We also defined a $\bdM$--equivariant injection 
\begin{equation}\label{eq:mapf}f\co H\otimes\lv{k+1}\to \frac{C_2(\lag{k+1})}{\lag{k+1}^{(1)}\wedge\lal{k+1}}.\end{equation}
We now recall this map.
We have an injection:
\[
\lv{k+1} \cong H_1(T_{k+1})\into H_1(T_{k+1};\R)\xrightarrow{v_*} H_1(\lal{k+1})\cong\lal{k+1}.
\]
Applying this injection in parallel with the one from equation~\eqref{eq:Hinj} gives us a map $H\otimes \lv{k+1}\into \lag{2}\otimes \lal{k+1}$.
Observe that $\lag{k+1}\wedge\lal{k+1}$ maps to $\lag{2}\otimes\lal{k+1}$ by $X\wedge Y\mapsto q(X)\otimes Y$, where $q\co \lag{k+1}\to\lag{2}$ is the projection induced by $\trunc{k+1}\to\trunc{2}$.
Since the kernel of this map is clearly $\lag{k+1}^{(1)}\wedge\lal{k+1}$, this gives us an isomorphism 
\[\lag{2}\otimes\lal{k+1}\cong\frac{\lag{k+1}\wedge\lal{k+1}}{\lag{k+1}^{(1)}\wedge\lal{k+1}}.\]
Define $f$ from equation~\eqref{eq:mapf} to be the composition of the map $H\otimes \lv{k+1}\to \lag{2}\otimes\lal{k+1}$ with the map above.
Recall that $\epsilon_k$ is the extension of the $k$th Morita homomorphism from Definition~\ref{de:emomap}. 
From Day~\cite[Main Theorem~B]{day}:
\begin{theorem}\label{th:dayb}
The map $\tilde d^2\circ \epsilon_k$ is a crossed homomorphism extending $f\circ \johom{k}$.
\end{theorem}

In this section we give a better interpretation of $\tilde d^2\circ \epsilon_k$.
Fix a map $j\co\Sigma\to X_{k+1}$ inducing the projection $\pi\to\trunc{k+1}$ on fundamental groups, as in Section~\ref{ss:nilstudy}.

\begin{proposition}\label{pr:surjh}
The crossed homomorphism $\tilde d^2\circ \epsilon_k$ is the principal crossed homomorphism of the element $[v(j_*[\Sigma])]\in C_2(\lag{k+1})/(\lag{k+1}^{(1)}\wedge\lal{k+1})$, where $[\Sigma]$ denotes any fundamental class of $\Sigma$ relative to $\partial\Sigma$.
\end{proposition}

\begin{proof}
Let $p\co X_{k+1}\to X_k$ be the projection and let $\rho\co\bdM\curvearrowright X_{k+1}$ denote the action from Proposition~\ref{pr:obs}.
For $[\phi]\in\bdM$, pick a representative $\phi\in \Diff(\Sigma,\partial\Sigma)$.
Pick a homotopy $F$ from $j$ to $\rho([\phi])\circ j\circ \phi^{-1}$ relative to $\partial\Sigma$.
Then $p\circ F$ is a homotopy from $p\circ j$ to $\rho([\phi])\circ p\circ j\circ \phi^{-1}$ relative to $\partial\Sigma$.
Let $T$ be a fundamental class of $\Sigma\times[0,1]$ relative to boundary.
By Definition~\ref{de:emomap}, $[v(p_*F_*T)]\in C_3(\lag{k})/B_3(\lag{k})$ is the value of the extended Morita map $\epsilon_k([\phi])$.

Note that $v(F_* T)\in C_3(\lag{k+1})$ is a lift of $v(p_*F_* T)$.
Let $i_0,i_1\co \Sigma\to\Sigma\times[0,1]$ be the endpoint inclusions.
We have fundamental classes $T_0$ and $T_1$ of $\Sigma$ relative to the boundary, and a chain $D\in C_2(\partial\Sigma\times[0,1])$, such that $\partial T= (i_1)_*T_1-(i_0)_*T_0+D$.
Since the restriction of $F$ to $\partial\Sigma\times[0,1]$ factors through $j|_{\partial\Sigma}$, Lemma~\ref{le:dim} implies that $v(F_*D)=0$ .
By Lemma~\ref{le:fc}, since $T_0$ and $\phi^{-1}_*T_1$ are both fundamental classes relative to the boundary, we have that $v(j_*T_0)=v(j_*\phi^{-1}_*T_1)$.
Putting this together, we have
\[
\begin{split}
\tilde d^2(\epsilon_k([\phi]))&=[\partial v(\tilde F_*T)]\\
&=[v(F_*(i_1)_*T_1)-v(F_*(i_0)_*T_0)+v(F_*D)]\\
&=[v(\rho(\phi)_*j_*\phi^{-1}_*T_1)-v(j_*T_0)]\\
&=[\phi]\cdot[v(j_*T_0)]-[v(j_*T_0)].\qedhere
\end{split}
\]
\end{proof}

Let $f'\co H\otimes \lv{k+1}\to C_2(\lag{k+1})/\bigwedge^2\lag{k+1}^{(1)}$ be the composition of $f$ above with the projection 
\[
\frac{C_2(\lag{k+1})}{(\lag{k+1}^{(1)}\wedge\lal{k+1})}\to \frac{C_2(\lag{k+1})}{\bigwedge^2\lag{k+1}^{(1)}}.
\]
\begin{corollary}\label{co:slightchange}
The homomorphism $f'$ is a  $\bdM$--equivariant injection and the principal crossed homomorphism $d[v(j_*[\Sigma])]$ of the element $[v(j_*[\Sigma])]\in C_2(\lag{k+1})/\bigwedge^2\lag{k+1}^{(1)}$ extends $f'\circ \johom{k}$.
\end{corollary}

\begin{proof}
The range of $f$ is $(\lag{k+1}\wedge\lal{k+1})/(\lag{k+1}^{(1)}\wedge\lal{k+1})$.
A computation using a basis for $\lag{k+1}$ shows that $(\lag{k+1}\wedge\lal{k+1})\cap \bigwedge^2\lag{k+1}^{(1)}$ is $\lag{k+1}^{(1)}\cap\lal{k+1}$.
This shows that $f'$ is injective.
The other statement follows immediately from Proposition~\ref{pr:surjh}.
\end{proof}

We proceed to frame this in a way similar to Theorem~\ref{th:mautfncobd}.
To do this, we will find a quotient of $C_2(\lag{k+1})$ in which $v(j_*[\Sigma])$ maps to a nontrivial, $\bdM$--invariant element.
The projection $q\co \lag{k}\to\lag{2}$ induces a projection $C_2(\lag{k})/\bigwedge^2\lag{k+1}^{(1)}\to C_2(\lag{2})$
(really the range is $C_2(\lag{2})/\bigwedge^2\lag{2}^{(1)}$ but $\lag{2}$ is abelian).
\begin{lemma}
The element $q_*v(j_*[\Sigma])\in C_2(\lag{2})$ is fixed by the action of $\bdM$.
\end{lemma}
\begin{proof}
Let $p\co X_{k+1}\to X_2$ be the projection induced from $\trunc{k+1}\to\trunc{2}$.
Then $q_*v(j_*[\Sigma])=v(p_*j_*[\Sigma])$.
By the argument in the proof Proposition~\ref{pr:surjh}, we know that for any $[\phi]\in\bdM$, the cycle 
$[\phi]\cdot v(p_*j_*[\Sigma])-v(p_*j_*[\Sigma])$ is a boundary in $B_2(\lag{2})$.
Since $\lag{2}$ is abelian, $B_2(\lag{2})=0$, proving the lemma.
\end{proof}

There is a slightly better way to characterize $q_*v(j_*[\Sigma])\in C_2(\lag{2})$.
\begin{definition}
The \emph{symplectic element} $\omega=\omega(\Sigma)$ of $\bigwedge^2 H\otimes \R$ is the image of the boundary element $[\partial\Sigma]\in[\pi,\pi]$ under the homomorphism $[a,b]\mapsto[a]\wedge[b]$.
If $a_1,\dotsc,a_g,b_1,\ldots,b_g$ are a basis for $\pi$ with $[\partial\Sigma]=[a_1,b_1]\dotsm[a_g,b_g]$, then
\[\omega=\sum_{i=1}^g[a_i]\wedge[b_i].\]
\end{definition}

\begin{lemma}\label{le:sigmatoomega}
There is a $\bdM$--equivariant isomorphism $C_2(\lag{2})\to\bigwedge^2H\otimes\R$ 
sending $q_*v(j_*[\Sigma])$ to the symplectic element $\omega$.
\end{lemma}

\begin{proof}
Since $\lag{2}$ is abelian, $C_2(\lag{2})$ and $H_2(\lag{2})$ are identical.
The inverse of $v_*$ maps $H_2(\lag{2})$ to $H_2(X_2;\R)$ isomorphically (by Corollary~\ref{co:vecisom}).
The map $\bigwedge^2 H^1(X_2;\R)\to H^2(X_2;\R)$ given by taking wedge products is an isomorphism since $X_2$ is a torus;
the evaluation pairings yield a dual isomorphism $H_2(X_2;\R)\to \bigwedge^2 H_1(X_2;\R)$.
Since $H_1(X_2)$ is torsion-free, we have $H_1(X_2;\R)\cong H_1(X_2)\otimes\R$.
The projection $\pi\to\trunc{2}$ induces an isomorphism $H=H_1(\pi)\cong H_1(\trunc{2})$, and we have an isomorphism $H_1(\trunc{2})\cong H_1(X_2)$ since $X_2$ is a $K(\trunc{2},1)$--space.
These isomorphisms are natural and therefore $\bdM$--equivariant.

Again let $p\co X_{k+1}\to X_2$ be the projection.
Since $p\circ j\co\Sigma\to X_2$ is a map inducing the canonical projection on fundamental groups and sending $\partial\Sigma$ to $[e]$, it descends to an Abel--Jacobi map from the closed surface $\Sigma/\partial\Sigma$ to $X_2$, which is a topological Jacobian torus for $\Sigma/\partial\Sigma$.
The push-forward of the fundamental class of $\Sigma/\partial\Sigma$ under an Abel--Jacobi map is the symplectic element in $H_2(X_2;\R)\cong \bigwedge^2H\otimes \R$ (this follows from considerations in Day~\cite[Section~2.3.2]{daysymp}). 
The integrals of left-invariant $2$--forms across the images of $\Sigma$ and $\Sigma/\partial\Sigma$ in $X_2$ are the same, so $q_*v(j_*[\Sigma])$ and $(p\circ j)_*[\Sigma/\partial\Sigma]$ map to the same element of $\bigwedge^2 H^1(X_2;\R)$.
\end{proof}

The kernel of $C_2(\lag{k+1})/\bigwedge^2\lag{k+1}^{(1)}\to \bigwedge^2 H\otimes\R$ is $(\lag{k+1}\wedge\lag{k+1}^{(1)})/\bigwedge^2\lag{k+1}^{(1)}$, which is isomorphic to $\lag{2}\otimes\lag{k+1}^{(1)}$.
By applying the $\bdM$--equivariant isomorphism $\lag{2}\cong H\otimes\R$ induced by $v$ to the first term, we have a $\bdM$--equivariant isomorphism $\lag{2}\otimes\lag{k+1}^{(1)}\to H\otimes\lag{k+1}^{(1)}$.
Thus the exact sequence:
\begin{equation}\label{eq:surfseq}0\to H\otimes \lag{k+1}^{(1)}\to \frac{C_2(\lag{k+1})}{\bigwedge^2\lag{k+1}^{(1)}}\to \bigwedge\nolimits^2 H\otimes\R\to 0.\end{equation}
Again, we consider the coboundary map from the resulting long exact sequence:
\[d\co H^0(\bdM;\bigwedge\nolimits^2 H\otimes\R)\to H^1(\bdM;H\otimes\lag{k+1}^{(1)}).\]
Let $\hat v\co H\otimes\lv{k+1}\to H\otimes \lag{k+1}^{(1)}$ be the $\bdM$--equivariant, injective homomorphism obtained by applying the injection from equation~\eqref{eq:Hinj} on the second part.
\begin{corollary}\label{co:surrange}
Any cocycle representing $d \omega\in H^1(\bdM;H\otimes\lag{k+1}^{(1)})$ is a crossed homomorphism extending $\hat v\circ \tau_k\co \bdI{k}\to H\otimes\lal{k+1}$.
\end{corollary}

\begin{proof}
By Lemma~\ref{le:sigmatoomega}, $\omega$ lifts to $[v(j_*[\Sigma])]\in C_2(\lag{k+1})/\bigwedge^2\lag{k+1}^{(1)}$.
By Corollary~\ref{co:slightchange}, the coboundary of $[v(j_*[\Sigma])]$ restricts to $f'\circ \johom{k}$.
Pulling back by the inclusion in sequence~\eqref{eq:surfseq}, this becomes $\hat v\circ \johom{k}$.
Any other representative differs from $d[v(j_*[\Sigma])]$ by $d\eta$, for some $\eta\in H\otimes\lag{k+1}^{(1)}$.
By Lemma~\ref{le:commact}, we know $\bdI{k}$ acts trivially on $H\otimes\lag{k+1}^{(1)}$, and $d\eta$ restricts to the zero map on~$\bdI{k}$.
\end{proof}

Theorem~\ref{th:mautfncobd} and Corollary~\ref{co:surrange} fit together nicely.
We map $C_2(\lag{k+1})=\bigwedge^2 \lag{k+1}$ to $\lag{2}\otimes\lag{k+1}$ by sending a basis element $X\wedge Y$ to $q(X)\otimes Y -q(Y)\otimes X$.
Since $\lag{2}\cong H\otimes \R$ (by tensoring the injection from equation~\eqref{eq:Hinj} with $\R$), Poincar\'e--Lefschetz duality induces a $\bdM$--equivariant isomorphism $\lag{2}\otimes\lag{k+1}$ to $\Hom(H,\lag{k+1})$.
The composition of these maps induces a $\bdM$--equivariant injection 
\[C_2(\lag{k+1})/(\bigwedge\nolimits^2\lag{k+1}^{(1)})\to \Hom(H,\lag{k+1}).\]
This injection induces a $\bdM$--equivariant injection 
\[\bigwedge\nolimits^2 H\otimes \R \to \Hom(H,H\otimes \R)\]
that sends the symplectic element $\omega$ to $\mathrm{id}\otimes 1$.
It also restricts to a $\bdM$--equivariant isomorphism between $H\otimes \lag{k+1}^{(1)}$ and $\Hom(H,\lag{k+1}^{(1)})$, giving us the following commutative diagram with exact rows:
\[
\xymatrix{
0 \ar[r] & H\otimes \lag{k+1}^{(1)} \ar[r]\ar[d]^{D_*} &
\frac{C_2(\lag{k+1})}{\bigwedge^2\lag{k+1}^{(1)}} \ar[r]\ar@{^{(}->}[d]^{D_*} & \bigwedge^2 H\otimes \R \ar[r]\ar@{^{(}->}[d]^{D_*} & 0\\
0 \ar[r] & \Hom(H,\lag{k+1}^{(1)})\ar[r] &\Hom(H,\lag{k+1})\ar[r]& \Hom(H,H\otimes\R) \ar[r] & 0\\
}
\]
Each of the vertical maps is induced by Poincar\'e--Lefschetz Duality; they are all injective and the first is an isomorphism.
Since this diagram is commutative and the vertical maps are $\bdM$--equivariant, we can relate our extension of $\hat v\circ \johom{k}\co \bdI{k}\to H\otimes \lal{k+1}$ to our extension of $\hat v\circ \fjh{k}\co \Aker{n}{k}\to \Hom(H,\lal{k+1})$.
The following needs no further proof:
\begin{corollary}\label{co:pd}
If $\tilde\omega\in C_2(\lag{k+1})/(\bigwedge^2\lag{k+1}^{(1)})$ is a lift of $\omega$, then $D_*\tilde\omega\in\Hom(H,\lag{k+1})$ is a lift of $\mathrm{id}\otimes1$, and the crossed homomorphisms $D_*\circ d\tilde\omega$ and $d(D_*\tilde\omega)$ are the same element of $Z^1(\bdM;\Hom(H,\lag{k+1}^{(1)}))$.
\end{corollary}

\subsection{Tightening the range}\label{ss:range}
Let $\latt_{k+1}$ be the $\Z$--span of $\log(\trunc{k+1})$ in $\lag{k+1}$.
This $\latt_{k+1}$ is an $(\Aut \pi)$--invariant (additive) subgroup, and is lattice in $\lag{k+1}$ by Theorem~\ref{th:lattice}.
In particular, it is finitely generated as an abelian group.

Let $S$ be a free generating set for $\pi$, and let $p\co \pi\to \trunc{k+1}$ be the projection.
Define $\jhob_S\co H\to \lag{k+1}$ by sending $[s]$ to $\log(p(s))$ for each $s\in S$, and extend linearly. 
Then by Lemma~\ref{le:vecabel}, this $\jhob_S$ is a lift of $\mathrm{id}\otimes 1$ whose image lies in $\latt_{k+1}$.
By Theorem~\ref{th:mautfncobd}, we have the following.
\begin{corollary}\label{co:nicest}
The crossed homomorphism 
\[d\jhob_S\in Z^1(\Aut \pi;\Hom(H,\lag{k+1}^{(1)}\cap\latt_{k+1}))\]
 extends $\hat v \circ \fjh{k}$, and its range is a finite-rank free abelian group.
\end{corollary}
Note that this proves Theorem~\ref{co:frab}, with $\eta=\hat v$, $A_k=\lag{k+1}^{(1)}\cap\latt_{k+1}$ and $\gamma_k=d\jhob_S$.

Now suppose that $n=2g$ and $S=\{a_1,\ldots,a_g$, $b_1,\ldots,b_g\}$ is a generating set for $\pi$, with $\prod_i[a_i,b_i]$ the class of the boundary loop in $\Sigma$.
Let $\tilde\omega_S$ be the image of $\sum_{i}\log(p(a_i))\wedge\log(p(b_i))$ in $C^2(\lag{k+1})/(\lag{k+1}^{(1)}\wedge\lal{k+1})$.
Let $D_*$ be as in Corollary~\ref{co:pd}.
\begin{corollary}\label{co:surnice}
We have that $D_*\tilde\omega_S=\jhob_S$.
Further, the crossed homomorphism 
\[\delta\tilde\omega_S\in Z^1(\bdM;D_*^{-1}(\Hom(H,\lag{k+1}^{(1)}\cap \latt_{k+1})))\]
 extends $\hat v\circ \fjh{k}$, and its range is a finite-rank free abelian group.
\end{corollary}

\begin{proof}
First we map $\tilde\omega_S$ to the element 
\[[\log(p(a_i))]\otimes\log(p(b_i)) -[\log(p(b_i))]\otimes \log(p(a_i)) \]
in $\lag{2}\otimes\lag{k+1}$.
By Lemma~\ref{le:vecabel}, the map $\log\circ p\co \pi\to \lag{2}$ induces the isomorphism $H_1(\pi)\cong H_1(\lag{2})=\lag{2}$.
Then Poincar\'e--Lefschetz duality on the first tensor component sends this element to $\jhob_S$.

Then by Corollary~\ref{co:pd}, we know that $D_*\circ \delta\tilde\omega_S$ is $\delta\jhob_S$.
By Corollary~\ref{co:nicest}, $D_*\circ\delta\tilde\omega_S$ has its range in $\Hom(H,\lag{k+1}^{(1)}\cap\latt_{k+1})$.
The range of $\delta\tilde\omega_S$ is the preimage $D_*^{-1}(\Hom(H,\lag{k+1}^{(1)}\cap \latt_{k+1}))$, which is a finite-rank free abelian group since $D_*$ is injective.
\end{proof}

\subsection{Example: the two-step case}\label{ss:example}
Let $\lag{}$ be a Lie algebra with its underlying vector space having a basis given by the $n+\binom{n}{2}$ vectors $\vct x_1,\dotsc,\vct x_n$, $\vct z_{12},\dotsc, \vct z_{n-1,n}$.
Let $\lag{}$ have the bracket operation with $[\vct x_i,\vct x_j]=\vct z_{ij}$ for $i<j$, and with $\vct z_{ij}$ central under the bracket.
Let $G$ be $\lag{}$ with the Baker--Campbell--Hausdorff product
\[\vct x * \vct y = \vct x+\vct y + \frac{1}{2}[\vct x,\vct y].\]
It is easy to show that $G$ is a contractible, two-step nilpotent Lie group.
For clarity, denote the element $\vct x_i$ by $\gamma_i$ when it is in $G$.
The logarithm map sends $\gamma_i\in G$ to $\vct x_i\in \lag{}$, and can be found on other elements using $\vct x * \vct y = \log (\exp(\vct x)\exp(\vct y))$.

Let $\Gamma=\trunc{3}$, with generators $a_1,\dotsc, a_n$.
The map $\Gamma\to G$ by $a_i\to \gamma_i$ extends to an injective homomorphism.
Let $\latt$ be the $\Z$--span of $\{\vct x_i\}_i\cup\{\frac{1}{2}\vct z_{ij}\}_{ij}$.
Then $\latt$ is a lattice in $\lag{}$ and contains $\log(\Gamma)$.
By Theorem~\ref{th:lattice}, $\Gamma$ is a lattice in $G$ and by Theorem~\ref{th:malcev}, $G$ is the Mal'cev completion of $\Gamma$.
Since $\log(a_ia_j)-\log(a_i)-\log(a_j)=\frac{1}{2}\vct z_{ij}$, we see that $\latt$ is the $\Z$--span of $\log(\Gamma)$.
Then in our earlier notation, $G=G_3$, $\lag{}=\lag{3}$ and $\latt=\latt_3$.

Let $S=\{\tilde a_i\}_i$ be generators for $\pi$ with $\tilde a_i$ mapping to $a_i\in\Gamma$.
Let $\xi_S$ be as in Corollary~\ref{co:nicest}.
\begin{example}\label{eg:transvect}
Let $\alpha\in\Aut\pi$ be the transvection $\tilde a_i\mapsto \tilde a_i\tilde a_j$ (other generators fixed).
Then $\alpha\cdot a_i = a_ia_j$, with other $a_l$ fixed.
We also see 
\[\begin{split}
\alpha \cdot \vct x_i &= \alpha\cdot \log (a_i) = \log (\alpha\cdot a_i) = \log (a_i a_j)\\
& = \log(\exp (\vct x_i)\exp(\vct x_j)) = \vct x_i*\vct x_j = \vct x_i+ \vct x_j +\frac{1}{2} \vct z_{ij}
\end{split}\]
 with other $\vct x_l$ fixed by $\alpha$.
To compute $(\alpha\cdot \xi_S) ([a_i])$, we send $[a_i]$ to $[a_i]-[a_j]$ (by $\alpha^{-1}\curvearrowright H$), then to $\vct x_i -\vct x_j$ by $\xi_S$, then to $\vct x_i+\vct x_j+\frac{1}{2}\vct z_{ij} - \vct x_j$ by $\alpha\curvearrowright \lag{}$.
So $(\alpha\cdot \xi_S)([a_i]) = \vct x_i +\frac{1}{2}\vct z_{ij}$ and $(\alpha\cdot \xi_S)([a_l])=\vct x_l$ for $l\neq i$.
Then since $\vct z_{ij}=[\log(a_i),\log(a_j)]$, we have:
\begin{equation}\label{eq:transvection}
d \xi_S(\alpha)([a_i]) = \frac{1}{2}[\log(a_i),\log(a_j)] \quad \text{and} \quad d \xi_S(\alpha)([a_l]) = 0
\end{equation}
for $l\neq i$.

If $\beta\in \Aut\pi$ is the inversion $\tilde a_i\mapsto \tilde a_i^{-1}$, with other generators fixed,
then a similar computation shows that $d \xi_S(\beta)=0$.

Transvections and inversions in $\{a_i\}_i$ generate $\Aut \pi$, so this is enough information to compute $d\xi_S$ on any element using the crossed homomorphism identity.
Note that $\hat v\co \Hom(H,\lv{3})\to \Hom(H,\lag{3}^{(1)})$ sends $[a_i]^*\otimes ([a_j,a_l])$ to $[a_i]^*\otimes [\log(a_j),\log(a_l)]$.
\end{example}

\begin{example}
Now suppose that $n=2g$.
Relabel $a_{g+1},\dotsc,a_{2g}$ as $b_1,\dotsc, b_g$, and similarly for $\{\tilde a_i\}_i$.
By a computation in Morita~\cite[Lemma~4.4]{mejo}, the images in $\Aut F_{2g}$ for the Lickorish generators for the mapping class group are given by two different classes of transvections, and the elements $\nu_i$ for $i=1,\dotsc, g-1$ as follows.
From Morita~\cite[Lemma~4.4]{mejo}:
\begin{gather*}
\nu_i(\tilde a_i) = \tilde a_i\tilde b_i^{-1}\tilde a_{i+1}\tilde b_{i+1}\tilde a_{i+1}^{-1}, \quad \nu_i(\tilde a_{i+1}) = \tilde a_{i+1}\tilde b_{i+1}^{-1}\tilde a_{i+1}^{-1}\tilde b_i\tilde a_{i+1}\\
\nu_i(\tilde b_i) = \tilde a_{i+1}\tilde b_{i+1}^{-1}\tilde a_{i+1}^{-1}\tilde b_i\tilde a_{i+1}\tilde b_{i+1}\tilde a_{i+1}^{-1},
\end{gather*}
where $\nu_i$ fixes the other generators.

As in Example~\ref{eg:transvect}, we can compute $d\xi_S(\nu_i)$ by computing the action of $\nu_i$ on the elements $\{\log(a_i)\}_i$ in $\lag{}$.
Again, this is easily accomplished using the Baker--Campbell--Hausdorff formula.
This computation, which we leave to the reader, shows that:
\begin{equation}\label{eq:Lickorishgen}
\begin{split}
d\xi_S(\nu_i)([a_i]) = &-\frac{1}{2}[\log(a_i),\log(b_i)]+[\log(a_{i+1}),\log(b_{i+1})]\\
 &+\frac{1}{2}[\log(a_i),\log(b_i)]+\frac{1}{2}[\log(b_i),\log(b_{i+1})]\\
d \xi_S(\nu_i)([a_{i+1}]) = & -\frac{1}{2}[\log(a_{i+1}),\log(b_{i+1})]\\
 & -\frac{1}{2}[\log(a_{i+1}),\log(b_{i})]-\frac{1}{2}[\log(b_{i}),\log(b_{i+1})]\\
d\xi_S(\nu_i)([b_i]) = & [\log(b_i),\log(b_{i+1})]
\end{split}
\end{equation}
and $d\xi_S(\nu_i)$ sends all other generators of $H$ to $0$.
There is a $\bdM$--equivariant injection $\bigwedge^2 H\into \latt_3\cap\lag{3}^{(1)}$ sending $[a_i]\wedge[a_j]\to [\log(a_i),\log(a_j)]$.
This induces a $\bdM$--equivariant injection $\Hom(H,\bigwedge^2 H)\into \Hom(H,\latt_3\cap\lag{3}^{(1)})$.
Comparing equations~\eqref{eq:transvection} and~\eqref{eq:Lickorishgen} to Lemmas~4.5 and~4.6 of Morita~\cite{mejo}, we see that this injection identifies $d\xi_S$ with Morita's crossed homomorphism~$\tilde k$.
\end{example}

\section{Extensions of Morita's homomorphisms}\label{se:moext}

\subsection{Polynomial chains for nilpotent homogeneous spaces}\label{ss:polynomialchains}
We start with some general theory that we apply in the next subsection.
Let $\Gamma$ be a finitely generated, torsion-free nilpotent group, $G$ its Mal'cev completion, $X=G/\Gamma$ the homogeneous space, and $\lag{}$ the Lie algebra of $G$.
The main result of this section is the construction of the \emph{polynomial straightening map}, an $(\Aut \Gamma)$--equivariant injection $s\co C_*(\Gamma)\to C_*(X)$ of the standard group homology chains $C_*(\Gamma)$ into the piecewise smooth singular homology chains $C_*(X)$.
This injection induces the canonical isomorphism $H_*(\Gamma)\cong H_*(X)$.
Roughly, to compute $s$ on a generator $C$ in $C_*(\Gamma)$, we pick a simplex $C'$ in $C_*(X)$ that realizes it geometrically, and straighten $C'$ to get a simplex that is parametrized by functions that are polynomials in exponential coordinates. 
Let $v\co C_*(X)\to C_*(\lag{})$ be the map Definition~\ref{de:vectorize}.
We show below that $v\circ s(C_*(\Gamma))$ is a finitely generated subcomplex of $C_*(\lag{})$.
Our result about extensions of Morita's homomorphisms follow easily.

As is standard, we model the $m$--simplex $\Delta_m$ as the subset of $\R^{m+1}$
\[\{(x_0,\ldots,x_m)|0\leq x_i\leq 1, 1= x_0+\ldots+x_m\}.\]
We label the vertices of $\Delta_m$ by setting $v_i$ to be the point with $x_i$--coordinate $1$ and all other coordinates $0$.
Let $A_i$ denote the face opposite $v_i$, \emph{i.e.} the subset of $\Delta_m$ on which $x_i=0$.

For $\gamma\in \Gamma$, recall the canonical representative $S^1\to X$ of $\gamma\in \pi_1(X)$ from Definition~\ref{de:canonical}.
It follows from the definitions that the lift $\R\to G$ of a canonical representative is a linear map in exponential coordinates.

\begin{definition}\label{de:polystr}
For $C\in C_m(\Gamma)$, the \emph{polynomial straightening} $s(C)\in C_m(X)$ is the map $s(C)\co \Delta_m\to X$ satisfying
\begin{itemize}
\item $s(C)$ sends the vertices of $\Delta_m$ to the basepoint $[e]$,
\item for each edge $E$ of $\Delta_m$, the loop $s(C)|_E$ is the canonical representative of the element of $\pi_1(X,[e])=\Gamma$ labeling $E$ in the $\Gamma$--labeling of $C$, and
\item for each $k, 1\leq k \leq m$, for each $k$--face $A$ of $\Delta_m$, $A\co \Delta_k\to\Delta_m$, the lift $\tilde {(s(C)\circ A)}\co\Delta_k\to G$ of $s(C)\circ A$ sending $v_0$ to the identity $e$ is a polynomial that is constant in $x_0$ and linear in $x_1$, in exponential coordinates.
\end{itemize}
\end{definition}

To show that polynomial straightenings exist, we wish to find a unique extension of polynomials defined on the faces of a simplex to the entire simplex.
The following lemma gives us these extensions.
\emph{A priori}, there is some ambiguity in such an extension. 
We address this with the second statement in the lemma.
This lemma is similar to part~(i) of the example on page~297 of Sullivan~\cite{infinitesimal}.

\begin{lemma}\label{le:polyext}
Suppose that $f_0,\ldots,f_m\co \R^{m+1}\to \R$ are polynomials 
such that $f_i$ agrees with $f_j$ on $A_i\cap A_j$ for all $i,j$.
Then there is a polynomial $f\co \R^{m+1}\to \R$ with 
$\deg f\leq \max\deg f_i$ such that $f$ agrees with $f_i$ on $A_i$ for each $i$.

Further, if $f_1,\ldots,f_m$ are all linear in $x_1$, then there is a unique such $f$ that is linear in $x_1$ and constant in $x_0$.
\end{lemma}

\begin{proof}
Let $\vct x$ denote $(x_0,\ldots,x_m)$, and for $k=0,\ldots, m$ let $p_k\co \R^{m+1}\to\R^{m+1}$ 
send $\vct x$ to $(x_0,\ldots,x_{k-1},0,x_{k+1},\ldots,x_m)$.

Since we only care about the values of the $f_i$ on $\Delta_m$, we may substitute $1-(\sum_{i=1}^mx_i)$ for $x_0$ and we may therefore assume that we are given $f_0,\ldots,f_m$ that are constant in $x_0$.
Define $\bar f_0=f_0$, and inductively assume we have defined $\bar f_{k-1}$ that agrees with $f_i$ on $A_i$ for $i=0,\ldots,k-1$ and is constant in $x_0$.
Then by the hypotheses, $f_k(\vct x)-\bar f_{k-1}(\vct x)$
evaluates to $0$ when $x_k=0$, $x_0=0$, and $\sum_{i=0}^mx_i=1$, \emph{i.e.} on $A_0\cap A_k$.

On $\langle A_k\rangle\cong \R^m$ (the subspace of $\R^{m+1}$ with $x_k=0$) the ideal of polynomials which vanish on $A_0\cap A_k$ is generated by the polynomials $x_0$ and $1-(x_1+\dotsb+x_{k-1}+x_{k+1}+\dotsb+x_m)$.
Since $f_k(p_k(\vct x))-\bar f_{k-1}(p_k(\vct x ))$ is constant in $x_0$ and $x_k$, we have that $1-(x_1+\dotsb+x_{k-1}+x_{k+1}+\dotsb+x_m)$ divides it.
So we may define $\bar f_k$ by:
\begin{equation}\label{eq:polyext}
\bar f_k(\vct x)=\bar f_{k-1}(\vct x)+(1-\sum_{i=1}^mx_i)\frac{f_k(p_k(\vct x))-\bar f_{k-1}(p_k(\vct x))}{1-(x_1+\ldots+x_{k-1}+x_{k+1}+\ldots+x_m)}.
\end{equation}
On $A_0$, the second term is zero, so $\bar f_k$ agrees with $\bar f_{k-1}$ and therefore with $f_0$.
On $A_k$, $x_k=0$ and the terms of $\bar f_k(\vct x)$ cancel each other out except for the $f_k(\vct x)$ term.
Therefore $\bar f_k$ agrees with $f_k$ on $A_k$.
For $A_i$, $1\leq i<k$, we have $x_i=0$, so $(x_0,\ldots,x_{k-1},0,x_{k+1},\ldots,x_m)$ is in $A_i\cap A_k$.
Then $f_k(x_0,\ldots,x_{k-1},0,x_{k+1},\ldots,x_m)=f_i(x_0,\ldots,x_{k-1},0,x_{k+1},\ldots,x_m)$, which equals $\bar f_{k-1}(x_0,\ldots,x_{k-1},0,x_{k+1},\ldots,x_m)$.
It follows that $\bar f_k$ agrees with $f_i$ on $A_i$.
This completes our inductive argument, and we take our final polynomial $f$ to be $\bar f_{m}$.

Now suppose that $f_1,\ldots,f_m$ are linear in $x_1$.
Since we only care about the value of $f_0$ on $A_0$, where $x_1=1-\sum_{i=2}^mx_i$, we may replace any instances of $x_1$ in $f_0(\vct x)$ with instances of $1-\sum_{i=2}^mx_i$.
So without loss of generality, we assume that $f_0$ is constant in $x_1$.
It then follows from the recursive construction of $f$ that each $\bar f_k$, $k=1,\ldots,m$, is linear in $x_1$.
Given our assumption that the $\{f_i\}$ are constant in $x_0$, it is also immediate that $f$ is constant in $x_0$.
If $\bar f$ is another polynomial that agrees with $f_i$ on $A_i$ for each $i$ and is constant in $x_0$, then $f-\bar f$ is in the ideal generated by $(1-\sum_{i=1}^mx_i)\prod_{i=1}^mx_i$.
Since any nonzero element of this ideal has degree $2$ or greater in $x_1$, the choice that is linear in $x_1$ is unique.
\end{proof}

\begin{proposition}\label{pr:strunique}
The polynomial straightening $s(C)$ exists and is unique for any $C\in C_m(\Gamma)$.
\end{proposition}

\begin{proof}
We proceed by induction on $m$.
Note that the polynomial straightening is the canonical representative if $m=1$, so it exists and is unique.

Now we suppose that we have unique polynomial straightenings of chains of dimension less than $m$.
Let $C_0,\ldots, C_m$ be the $(m-1)$--faces of $C$.
For each $i$, we get a straightening $s(C_i)\co \Delta_{m-1}\to X$.
For $i=1,\ldots,m$, let $\tilde s(C_i)\co \Delta_{m-1}\to G$ be the lift sending $v_0$ to the identity; let $\tilde s(C_0)\co\Delta_{m-1}\to G$ be the lift sending $v_0$ to $\tilde s(C_1)(v_1)\in G$.
Let $p_i\co \Delta_{m}\to\Delta_{m-1}$ be the linear projection (degeneracy map) to the $i$th face.
For each $i$, set $f_i=s(C_i)\circ p_i\co \Delta_m\to G$.

If $i>0$, then it follows from the consistency of the labels on $C$ that $f_i$ and $f_0$ both restrict to the same lift of the straightening of the same $(m-2)$ face of $C$ on $A_0\cap A_i$.
Since polynomial straightenings for $(m-2)$--chains are unique, we have that $f_i$ agrees with $f_j$ on $A_i\cap A_j$ for all $i, j$.
Then we apply Lemma~\ref{le:polyext} (coordinate by coordinate) to get a map $f\co \Delta_m\to G$, polynomial in exponential coordinates, that agrees with $f_i$ on $A_i$ for all $i$ and is constant in $x_0$ and linear in $x_1$.
This $f$ descends to our polynomial straightening $s(C)\co \Delta_m\to X$.
The properties in the definition follow by construction.

The uniqueness of $s(C)$ follows from the uniqueness in Lemma~\ref{le:polyext} (using the third property in the definition) and inductively from the uniqueness of the straightenings of the faces.
\end{proof}

\begin{proposition}\label{pr:strequ}
The polynomial straightening map $s\co C_*(\Gamma)\to C_*(X)$ is a chain map inducing the canonical isomorphism $H_*(\Gamma)\cong H_*(X)$ and is $(\Aut \Gamma)$--equivariant.
\end{proposition}

\begin{proof}
Let $C_i$ be the $i$th face of $C\in C_*(\Gamma)$ and let $D_i$ be the $i$th face of $s(C)$.
Then by definition, $s(C_i)$ and $D_i$ are both polynomial straightenings of $C_i$, so by Proposition~\ref{pr:strunique}, $s(C_i)=D_i$.
It follows immediately that $s$ intertwines the boundary operations of $C_*(\Gamma)$ and $C_*(X)$.

Let $F_*(\Gamma)$ be the standard resolution of $\Z$ over $\Z\Gamma$; $F_m(\Gamma)$ is the free $\Gamma$--module generated by the set of $\Gamma$--labeled $m$--simplices.
The singular chains $C_*(G)$ are also a resolution of $\Z$ over $\Z\Gamma$.
For $\gamma\cdot C\in F_m(\Gamma)$, let $\tilde s(\gamma\cdot C)\co \Delta_m\to G$ be the lift of $s(C)\to \Delta_m\to X$ that sends $v_0$ to $\gamma$.
It follows from standard arguments (see Brown~\cite[Chapter~I.7]{brown}) that $\tilde s\co F_*(\Gamma)\to C_*(G)$ is a homotopy equivalence of $\Gamma$--complexes.
Therefore its induced map on chains $s\co C_*(\Gamma)\to C_*(X)$ in turn induces an isomorphism on homology.

The $(\Aut \Gamma)$--equivariance follows from Proposition~\ref{pr:strunique}: for $\phi\in\Aut\Gamma$ and $C\in C_*(\Gamma)$, the chain $\phi\cdot s(C)$ is a polynomial straightening of $\phi\cdot C$ and by uniqueness must equal $s(\phi\cdot C)$.
\end{proof}

Let $\latt$ be the $\Z$--span of $\log(\Gamma)$.
Let $\tilde s(C)\co \R^{k+1}\to G$ be the lift of $s(C)\co \R^{k+1}\to X$ sending $v_0\in\Delta_k$ to~$e$.
\begin{lemma}\label{le:KfortildeC}
There is an integer $K>0$ depending only on $\Gamma$ and $k$ such that 
for any basis $B$ for $\latt$ and any generator $C\in C_k(\Gamma)$, the coordinate polynomials for $\log\circ \tilde s(C)\co \R^{k+1}\to \lag{}$ with respect to $B$ have coefficients in~$\frac{1}{K}\Z$.
\end{lemma}

\begin{proof}
We induct on $k$.
Let $C\in C_1(\Gamma)$ be a $\Gamma$--labeled $1$--simplex.
Then $\log\circ \tilde s(C)\co \R^2\to \lag{}$ is constant in the first coordinate and is given by a $1$-parameter family map in the second coordinate.
In exponential coordinates, $1$-parameter family maps are linear; since this is a linear map sending $v_1\in\Z^2$ to an element of $\latt$, this map has integer coefficients in its coordinate polynomials with respect to any basis for $\latt$.
This proves the case $k=1$.

Suppose that the lemma holds for $k-1$ with a constant $K_0>0$.
Fix a basis $B$ for $\latt$.
Suppose $C\in C_k(\Gamma)$ is a generator.
For $i=1,\dotsc,k$, the restriction of $\tilde s(C)$ to $\langle A_i\rangle$ (the span of $A_i$) is $\tilde s(C_i)$, where $C_i$ is the $i$th face of $C$.
Therefore the coordinate polynomials for these restrictions with respect to $B$ have coefficients in $\frac{1}{K_0}\Z$.
Let $f_0$ denote the restriction of $\tilde s(C)$ to $\langle A_0\rangle$.
Then $f_0$ is $\lmult{g}\circ \tilde s(C_0)$, where $C_0$ is the $0$th face of $C$ and $g=\tilde s(C)(v_1)$.
By Corollaries~\ref{co:polymult} and~\ref{co:multlattice}, and Proposition~\ref{pr:integervaluedcoeffs}, $\lmult{g}$ is a polynomial in exponential coordinates, and its coefficients with respect to $B$ are in $\frac{1}{K_1}\Z$ for some integer $K_1>0$ depending only on $\Gamma$.
Further, the polynomial degree of $\lmult{g}$ is bounded by the nilpotence class of $\Gamma$.
Then by Lemma~\ref{le:polycomp} there is an integer $K>0$, depending only on $k$ and $\Gamma$, such that the coefficients of $f_0$ with respect to $B$ are in $\frac{1}{K}\Z$.
Then by equation~\eqref{eq:polyext}, the coefficients of $\tilde s(C)$ with respect to $B$ are no worse than those of $f_0$, and are also in $\frac{1}{K}\Z$.
\end{proof}

\begin{proposition}\label{pr:Kexists}
Let $k$ be an integer.
There is a positive integer $K$ depending only on $k$ and $\Gamma$, such that for any
polynomially straightened $k$--simplex $C\co \Delta^k\to X$, the value
$v(C)$ is in the lattice $\frac{1}{K}\bigwedge^k\latt$ in $C_k(\lag{})$.
\end{proposition}

\begin{proof}
For $g\in G$, let $\lmult{g}\co G\to G$ denote left-multiplication by $G$.
Let $\alpha$ be a left-invariant $k$--form on $G$.
For $g\in G$ and a tangent multi-vector $\vct{x}\in \bigwedge^k T_g^*G$, we can express the value of $\alpha$ at $g$ on $\vct{x}$ in terms of $\alpha$ at $e$ and the derivative of the left-multiplication function:
\[\alpha_g(\vct{x})=\alpha_e(D_g\lmult{g^{-1}}(\vct{x})).\]
Let $\tilde C\co \Delta^k\to G$ be the lift of $C$ sending $v_0$ to $e$.
Then for $p\in\Delta^k$ and $\vct{x}\in \bigwedge^k T_p\Delta^k$, we have
\begin{equation}\label{eq:polypullback}
(\tilde C^*\alpha)_p(\vct{x})=\alpha_e(D_{\tilde C(p)}\lmult{\tilde C(p)^{-1}}(D_p\tilde C(\vct{x}))).
\end{equation}

Let $B$ be a basis for $\bigwedge^k\latt$ and let $B^*$ be its dual basis in $\bigwedge^k\latt^*$.
Now we suppose that $\alpha$ is a left-invariant $k$--form on $G$ such that $\alpha_e\in B^*$.
Then $\alpha_e\co \bigwedge^k\lag{}\to \R$ is a linear map sending $\bigwedge^k\latt$ to $\Z$.

Consider the map $f\co G\times \bigwedge^kTG\to \bigwedge^kTG$ given by sending the point $g$ and the tangent $k$--vector $\vct{x}$  at the point $h$ to the tangent $k$--vector $D_h\lmult{g}(\vct{x})$ at the point $g h$.
By Corollary~\ref{co:polymult} and Lemma~\ref{le:polyderiv}, this map is a polynomial map in exponential coordinates.
Further, the polynomial degree of this map is bounded by the nilpotence class of $G$.
Since $\bigwedge^k T\lag{}=\lag{}\times\bigwedge^k\lag{}$, we may think of $f$ in exponential coordinates as a map $f\co \lag{}\times\lag{}\times\bigwedge^k \lag{}\to\lag{}\times\bigwedge^k\lag{}$.
Again by Corollary~\ref{co:polymult} and Lemma~\ref{le:polyderiv}, there is an integer $K_0>0$ such that $f$ sends $\latt\times\latt\times\bigwedge^k\latt$ into $\latt\times\frac{1}{K_0}\bigwedge^k\latt$.
This $K_0$ depends only on $k$ and $\Gamma$.

By definition, the map $\tilde C\co \Delta^k\to G$ extends to a map $\tilde C\co \R^{k+1}\to G$ that is polynomial in exponential coordinates.
By Lemma~\ref{le:KfortildeC}, there is an integer $K_1>0$ such that $\tilde C$ sends $\Z^{k+1}$ to $\frac{1}{K_1}\latt$ in exponential coordinates.
This $K_1$ depends only on $k$ and $\Gamma$.
By Lemma~\ref{le:polyderiv}, there is an integer $K_2>0$, depending only on $k$ and $\Gamma$, such that $D\tilde C\co \bigwedge^kT\R^k\to\bigwedge^kTG$ sends $\Z^{k+1}\times\bigwedge^k\Z^{k+1}$ to $\frac{1}{K_1}\latt\times\frac{1}{K_2}\bigwedge^k\latt$ in exponential coordinates.

Finally, consider the map $h\co \bigwedge^kTG\to G\times \bigwedge^kTG$ that sends the tangent $k$--vector $\vct{x}$ at the point $g$ to $(g^{-1},\vct{x})$.
In exponential coordinates, this is a linear map $\lag{}\times\lag{}\times\bigwedge^k\lag{}\to\lag{}\times\bigwedge^k\lag{}$ sending $\latt\times\latt\times\bigwedge^k\latt$ to $\latt\times\bigwedge^k\latt$.

By equation~\eqref{eq:polypullback}, the map $(p,\vct{x})\mapsto (\tilde C^*\alpha)_p(\vct{x})$ is the composition $\alpha_e\circ f\circ h\circ D \tilde C$.
By Lemma~\ref{le:polycomp}, there is an integer $K_3>0$ such that for any $p,\vct{x}\in \Z^k$, we have $\tilde C^*\alpha_p(\vct{x})\in\frac{1}{K_3}\Z$.
This $K_3$ ultimately only depends on $k$ and $\Gamma$.
By Proposition~\ref{pr:integervaluedcoeffs}, there is an integer $K_4$ depending only on $k$ and $\Gamma$ such that $\tilde C^*\alpha\co T\R^k\to \R$ is a polynomial with coefficients in $\frac{1}{K_4}\Z$.

We use Proposition~\ref{pr:basisform} to evaluate $v(C)$.
To find the coefficients of $v(C)$ in the basis expression with respect to $B$, we integrate $\alpha^*\tilde C$ over the standard $k$--simplex $\Delta^k$ for $\alpha$ with $\alpha_e$ ranging over $B^*$.
We compute the integral as an iterated integral on $\alpha^*\tilde C$.
\[\int_{\tilde C}\alpha=\int_0^1\int_0^{x_1}\cdots\int_0^{x_{m-1}} (\tilde C^*\alpha)_{(x_1,\ldots,x_m)}(\frac{\partial}{\partial x_1},\ldots,\frac{\partial}{\partial x_m})dx_m\cdots dx_1.\]
Integrating with respect to each variable, we divide the coefficients by the appropriate degree and evaluate.
The end result is that the value of the integral is in $\frac{1}{K}\Z$, where the integer $K>0$ depends only on $K_4$ and $k$.
Then by Proposition~\ref{pr:basisform}, we see that $v(C)\in\frac{1}{K}\bigwedge^k\latt$.
\end{proof}

\subsection{A more careful extension of Morita's homomorphisms}\label{ss:mocareful}
Now we return to the situation where $\Gamma=\trunc{k}$ is the class--$(k-1)$ nilpotent truncation of the fundamental group $\pi$ of the once-bounded, genus--$g$ surface $\Sigma$.
The goal of this section is to prove Theorem~\ref{th:mobetter}.
Let $p\co \pi\to\trunc{k}$ denote the projection.

Let $K>0$ be the smallest integer such that the image of $v\circ s\co C_*(\trunc{k})\to C_*(\lag{})$ is contained in $\frac{1}{K}\bigwedge^3\latt_k$.
This $K$ exists by Proposition~\ref{pr:Kexists}.
Let $B_k$ denote the $\bdM$--module
\[B_k=\frac{\frac{1}{K}\bigwedge\nolimits^3\latt_k+B_3(\lag{k})}{B_3(\lag{k})}.\]
\begin{definition} \label{de:mobetter}
Let $[\ell]\in C_1(\pi)$ and $C_\Sigma\in C_2(\pi)$ be as in Definition~\ref{de:mobetter}.
For $[\phi]\in\bdM$, let $D_{[\phi]}\in C_3(\pi)$ with $\partial D_{[\phi]}=[\phi]_*C_\Sigma-C_\Sigma$.  
The \emph{polynomial-straightening extended Morita map}  is 
\[\begin{split}
\psmo\co \bdM&\to B_k\\
[\phi]&\mapsto v(s(p_*D_{[\phi]})).\end{split}\]
\end{definition}

\begin{lemma}\label{le:choices}
If $C_\Sigma$ and $C_\Sigma'\in C_2(\pi)$ are as in Definition~\ref{de:mobetter},
then there is a chain $D\in C_3(\pi)$ with 
\begin{equation}
\label{eq:bddiff}
\partial D = C_\Sigma-C_\Sigma'.
\end{equation}
For any $[\phi]\in\bdM$, let $D_{[\phi]}$ and $D_{[\phi]}'\in C_3(\pi)$ with $\partial D_{[\phi]}=[\phi]_*C_\Sigma-C_\Sigma$ and $\partial D_{[\phi]}'=[\phi]_*C_\Sigma'-C_\Sigma'$.
Then 
\begin{equation}
\label{eq:cobddiff}
D_{[\phi]}-D_{[\phi]}'-([\phi]_*D-D)\in B_3(\pi).
\end{equation}
If $C_\Sigma=C_\Sigma'$ then we can take $D=0$.
\end{lemma}

\begin{proof}
Recall that $H_m(\pi)=0$ for $m>1$.
The chain $C_\Sigma-C_\Sigma'$ is in $Z_2(\pi)$ and therefore in $B_2(\pi)$.
Pick any $D\in C_3(\pi)$ satisfying equation~\eqref{eq:bddiff}.
It is immediate that the chain in equation~\eqref{eq:cobddiff} is in $Z_3(\pi)$, and therefore in $B_3(\pi)$.
\end{proof}

\begin{proposition}\label{pr:mocarefulwelldef}
For a fixed $C_\Sigma\in C_2(\pi)$, the map $\psmo$ is well defined.
\end{proposition}
\begin{proof}
Let $D_{[\phi]}$ and $D_{[\phi]}'$ be as in Lemma~\ref{le:choices}, with $C_\Sigma=C_\Sigma'$.
Then by that lemma, $D_{[\phi]}-D_{[\phi]}'\in B_3(\pi)$.
Since $v\circ s\circ p_*$ is a chain map, we have that $v(s(p_*D_{[\phi]}))-v(s(p_*D_{[\phi]}'))\in B_3(\lag{k})$.
\end{proof}

\begin{proposition}
The cohomology class
\[[\psmo]\in H^1(\bdM;B_k)\]
does not depend on any of the choices.
\end{proposition}

\begin{proof}
Let $C_\Sigma$, $C_\Sigma'$, $D$, $[\phi]$, $D_{[\phi]}$, and $D_{[\phi]}'$ be as in Lemma~\ref{le:choices}.
Then since $v\circ s\circ p_*$ is a chain map, equation~\eqref{eq:cobddiff} implies that for any $[\phi]$
\[v(s(p_*D_{[\phi]}))-v(s(p_*D_{[\phi]}'))-\big([\phi]\cdot v(s(p_*D))-v(s(p_*D))\big)\in B_3(\lag{k}).\]
In particular, the first term represents $\psmo$ computed with $C_\Sigma$, the second term represents $\psmo$ computed with $C_\Sigma'$, and the third term represents the coboundary $d [v(s(p_*D))]$ evaluated on $[\phi]$.
\end{proof}

\begin{proposition}\label{pr:lattice}
The $\bdM$--module $B_k$
is a finitely generated abelian group and a lattice in $C_3(\lag{k})/B_3(\lag{k})$.
\end{proposition}

\begin{proof}
The module $B_k$ is isomorphic to
\[\frac{\frac{1}{K}\bigwedge\nolimits^3\latt_k}{(\frac{1}{K}\bigwedge\nolimits^3\latt_k)\cap B_3(\lag{k})}\]
and is therefore finitely generated.

By the formula for $\partial\co C_m(\lag{k})\to C_{m-1}(\lag{k})$ and the fact that $\lag{k}$ has rational structure constants with respect to any basis for $\latt$, we see that 
$\partial$ restricts to $\partial\co \bigwedge^m\latt_k\otimes\Q\to \bigwedge^{m-1}\latt_{k}\otimes\Q$.
Therefore the splitting $C_3(\lag{k})\cong B_3(\lag{k})\oplus (C_3(\lag{k})/B_3(\lag{k}))$ is rational with respect to $\frac{1}{K}\bigwedge^3\latt_k$.
In particular, since $\frac{1}{K}\bigwedge\latt_k$ is a lattice in $C_3(\lag{k})$, its image is a lattice in $C_3(\lag{k})/B_3(\lag{k})$.
\end{proof}

\begin{proof}[Proof of Theorem~\ref{th:mobetter}]
Let $\eta\co H_3(\Gamma)\to B_k$ be the composition 
\[H_3(\Gamma)\xrightarrow{s_*} H_3(X)\xrightarrow{v_*} B_k.\]
Since $s$ and $v$ are $\bdM$--equivariant, $\eta$ is as well.
We need to show that $\eta$ is injective and that $\psmo$ extends $\eta\circ \mohom{k}$.

By Igusa--Orr~\cite[Theorem~5.9]{igusaorr}, $H_3(\Gamma)$ is torsion-free.
Therefore $H_3(\Gamma)\to H_3(X;\R)$ is injective.
The map $H_3(X)\xrightarrow{v_*}B_k\cap H_3(\lag{k})$ is a factor of the injection $H_3(X)\to H_3(X;\R)\xrightarrow{v_*}H_3(\lag{k})$ and is therefore injective.
So $\eta$ is an injection.

The association $\bdM\to C_3(\pi)$ given by $[\phi]\mapsto D_{[\phi]}$ is the same in Definition~\ref{de:mohom} and Definition~\ref{de:mobetter}.
From this it follows immediately that $\psmo$ extends~$\eta\circ\mohom{k}$.
\end{proof}

\bibliographystyle{amsplain}
\bibliography{autfnejh}

\noindent
Department of Mathematics 253-37\\
California Institute of Technology\\
Pasadena, CA 91125\\
E-mail: {\tt mattday@caltech.edu}
\medskip

\end{document}